\documentclass[preprint,12pt]{elsarticle}
\usepackage{comment}
\usepackage{graphicx}
\usepackage{subfigure}
\usepackage{amsmath,amsfonts,amssymb,amsbsy,stmaryrd}
\usepackage{lineno}
\usepackage{psfrag}
\graphicspath{{figures/}}

\newcommand\hooke{\mathbb{H}}
\newcommand\trace{\ensuremath{\operatorname{tr}}}
\newcommand\KA[1]{{\mathrm{KA}}\!\left(#1\right)}

\newcommand\KAoo[1]{{\mathrm{KA}}^{00}\!\left(#1\right)}
\newcommand\KAh[1]{\mathrm{KA}_h\!\left(#1\right)}

\newcommand\KAhoo[1]{\mathrm{KA}_h^{00}\!\left(#1\right)}
\newcommand\SA[1]{\mathop{\mathrm{SA}}\!\left(#1\right)}

\newcommand\ecr[2]{\mathop{\mathrm{e_{CR(#2)}}}\!\left(#1\right)}

\newcommand\enernorm[2]{\|#1\|_{\hooke^{-1},#2}}
\newcommand\strainnorm[2]{\|#1\|_{\hooke,#2}}

\newcommand{\globalerror}{\ensuremath{e_h}}

\newcommand{\s}{\ensuremath{^{(s)}}}
\newcommand{\g}{\ensuremath{^\square}}
\newcommand\assemg{\mathcal{A}}

\newcommand\shapef{\varphi}
\newcommand\shapev{\boldsymbol{\shapef}_h}

\newcommand\admiss[1]{\widehat{#1}}

\newcommand{\dep}{\ensuremath{\mathbf{u}}}
\newcommand{\hdep}{\ensuremath{\admiss{\dep}}}
\newcommand{\stiff}{\ensuremath{\mathbf{K}}}
\newcommand{\force}{\ensuremath{\mathbf{f}}}
\newcommand{\F}{\ensuremath{\mathbf{F}}}
\newcommand{\hF}{\ensuremath{\admiss{\F}}}

\newcommand{\adF}{\ensuremath{\admiss{F}}}

\newcommand{\lam}{\ensuremath{\boldsymbol{\lambda}}}

\newcommand{\hlam}{\ensuremath{\admiss{\lam}}}

\newcommand{\traceh}{\ensuremath{\mathbf{t}}}

\newcommand{\passem}{\ensuremath{\mathbf{A}}}
\newcommand{\dassem}{\ensuremath{\mathbf{\underline{A}}}}

\newcommand{\schur}{\ensuremath{\mathbf{S}}}

\newcommand{\rhs}{\ensuremath{\mathbf{b}}}
\newcommand{\kernel}{\ensuremath{\mathbf{R}}}

\newcommand{\alp}{\ensuremath{\boldsymbol{\alpha}}}

\newcommand{\psassem}{\ensuremath{\tilde{\passem}}}
\newcommand{\dsassem}{\ensuremath{\tilde{\dassem}}}

\newcommand\strain[1]{\varepsilon\left(#1\right)}

\newcommand\ecrseq{\mathrm{e_{CR}^{seq}}}
\newcommand\ecrpara{\mathrm{e_{CR}^{ddm}}}
\newcommand\ecrparaE{\mathrm{e_{CR}^{ddm,E}}}

\newcommand\structure{\Omega}

\title{Fast estimation of discretization error for FE problems solved by domain decomposition}
\author[lmt]{A.˜Parret-Fr\'eaud}
\ead{augustin.parret-freaud@lmt.ens-cachan.fr}
\author[lmt]{C.~Rey\corref{cor1}}
\ead{christian.rey@lmt.ens-cachan.fr}
\author[lmt]{P.~Gosselet}
\ead{pierre.gosselet@lmt.ens-cachan.fr}
\author[onera]{F.~Feyel}
\ead{frederic.feyel@onera.fr}
\cortext[cor1]{Corresponding author}
\address[lmt]{LMT-Cachan, ENS Cachan/CNRS/UPMC/PRES UniverSud, 61 av. du pr\'esident Wilson, 94235 Cachan cedex, France}
\address[onera]{ONERA, DMSM/CEMN, 29 avenue de la division Leclerc, BP72, F92322 Chatillon cedex, France}

\begin{document}

\begin{abstract}
  This  paper  presents  a strategy  for a posteriori error estimation  for 
  substructured  problems solved  by  non-overlapping  domain decomposition methods. 
  We focus on global estimates of the discretization error obtained through the error  in constitutive relation for linear  mechanical problems.   
  Our method allows to  compute error estimate in a  fully parallel way for  both primal (BDD) and 
  dual (FETI) approaches of  non-overlapping domain  decomposition whatever the state
  (converged or not) of the associated iterative solver. 
Results obtained on an academic problem show that the strategy we propose is efficient in the sense that correct estimation is obtained with fully parallel computations; they also indicate that the estimation of the discretization error reaches sufficient precision in very few iterations of the domain decomposition solver, which enables to consider highly effective adaptive computational strategies.
%One may then be able to
  %access to the error estimate at each iteration of the global computation. When
  %solver's convergence  is reached, our  method gives good results  compared to
  %the  estimate computed  through a  sequential  problem.  Moreover,  it can  be
  %observed that convergence of global error  estimate is far faster than the one
  %related  to the solver.  This last  property is  also observed  on the  map of
  %elementary  contributions  to  the  global  error.
  %,  allowing  to  use  it  for adaptivity  purpose.   
  %The  numerical  performance  is  illustrated  through results on one academic problem.
\end{abstract}

\maketitle

\section{Introduction}
\label{sec:intro}

The  setting-up   of  robust  numerical   methods  to  solve   complex systems of partial
differential  equations  has become  a  key  issue  in applied  mathematics  and
engineering,  driven by  the  increasing  use of  numerical  simulation in  both
research and industry. Among the latter, virtual testing has become a short term
aim, with the objective to replace expensive experimental  studies and  validations by
numerical simulations, even  in order to certify large  structures as planes and
bridges.

Thus,  one  key   point  of  the   numerical  methods   to  develop  is  the
\emph{verification} of computations which  enables to warranty that the computed
solution is sufficiently close to  the original continuum mechanics model.  This
topic of  numerical analysis has been the  subject of many studies  for the last
decades.   Three main  classes of  error  estimator have  been developed,  based
either  on   equilibrium  residuals  \cite{babuska:78:resEq},   flux  projection
\cite{zienkiewicz:87:zz1}      or       error      in      constitutive      law
\cite{ladeveze:75:erreur}.  An overview of those various methods can be found in
\cite{ladeveze:04:mastCalculations}.

Another key point  of numerical methods is their  ability to quickly provide
solutions to large (nonlinear) systems.  The most classical answer to this issue
is  to use  domain  decomposition methods  in  order to  take  advantage of  the
parallel hardware  architecture of recent  clusters and grids.   In engineering,
non-overlapping domain  decomposition methods are  mostly employed, such  as the
well known FETI \cite{FARHAT:1994:ADV} or BDD \cite{mandel:93:ddm}.  An overview
of the  main approaches related  to non-overlapping domain decomposition  can be
found in \cite{GOSSELET.2007.1}.

We  aim  to provide  fully  integrated  adaptive strategies  to
compute large structural mechanics problems  with certified quality. To do that,
our  current  approach   is  to  explore  some  ways   of  making  bidirectional
interactions  between domain  decomposition and  a posteriori  error estimation.
Our developments are based both on the error in constitutive relation to measure
the quality  of our results  and to forecast  mesh refinement, and on  a generic
vision  of non-overlapping  domain  decomposition methods  which  enables to  do
high-performance computing.

This  paper  focuses on  the  estimation of  the  global  error in  constitutive relation in order (among  others) to study how it is influenced  by the error in the  convergence of  the  domain decomposition  solver  which is  linked to the non-satisfaction of interface equations (continuity of displacements and balance of forces). To do so we  propose a strategy to  build, in parallel  and %in the course of 
during the iterations, displacement and stress fields which are kinematically admissible (KA) and  statically   admissible (SA) on  the   whole  structure.   We  face   two  main difficulties.  First,  since before convergence interface fields  do not possess the classical properties of  discretized fields (continuity of displacements and weak  equilibrium),  the  recovery  of  admissible  displacements  and  stresses requires some preprocessing.  Second, the  computation of statically admissible fields being an operation which can not  be conducted independently on each element (in some methods  it can  even be a large bandwidth operation),  classical recovery methods      \cite{ladeveze:leguillon:81:erdc:sarecovery, ladeveze:rougeot:1997:erdc:sarecovery, pares:diez:huerta:2006:fluxfreesubdomain, gallimard:2008:erdc:sarecovery, mointinhoDeAlmeida:maunder:2009:sarecoverypu, ladeveze:chamoin:florentin:2009:erdc:sarecovery} would require inter-subdomain communications. 

Our generic method to build continuous displacement and balanced traction fields for both  primal and  dual  approaches of  non-overlapping  domain decomposition  is presented  through this  paper.  It  will be  shown that  the properties  of the preconditioners   involved   in   domain   decomposition   solvers   make   this reconstruction costless, and  that an error estimator can then  be computed in a fully parallel way.

%The scope of this paper  is the following.  
This paper is organized as follows. Section \ref{sec:basics} recalls the general framework related to  our upcoming developments, mainly the estimation of the error in constitutive equation and the use of domain decomposition method. %A classical mechanical
% problem is formulated together with  its substructured approach, in a algebraic
% way that directly  use admissibility spaces and error  in constitutive relation.
% Then, basic  principles of a  posteriori error estimation based  on constitutive
% relation is  presented after having introduced the  finite element approximation
% of the  previous problem. At  last, the main  steps of the SA stress recovery
% process proposed in \cite{ladeveze:leguillon:81:erdc:sarecovery} are quickly  recalled. 
Section \ref{sec:errordd} shows how the problem of error estimation in a substructured context can be brought back to the computation of nodal displacement and traction fields which are admissible in a discrete sense.
Sections  \ref{sec:bddrecovery}  and \ref{sec:fetirecovery}  describes  how to obtain these fields without inter-subdomains exchanges when using classical primal (BDD) and dual (FETI) domain decomposition methods with good preconditioners.
% to  build  continuous displacement  and  balanced  stress  fields  in  a  parallel  context,  focusing successively on  the specificity of  the primal (\ref{sec:bddrecovery})  and the
%dual (\ref{sec:fetirecovery}) approach  of non-overlapping domain decomposition.
%Section  \ref{sec:sarecovery}  presents  a  simple  way  to  compute  continuous balanced  traction field  on  interfaces used  further  to recover  equilibrated stress  field   on  sub-domains.   
Section   \ref{sec:numeric}  presents numerical assessments, first to validate the parallel recovery procedure, then to prove that a good estimation can be obtained far earlier than the solver converged (in the sense of domain decomposition iterative solver).
%then on the  convergence of error estimate  towards the iterations of the solver.  
Finally, Section \ref{sec:conclusions} concludes this paper.

\section{Framework of the study}
\label{sec:basics}

\subsection{Reference mechanical problem}
\label{sec:crebasicsrefpb}

Let us consider the static equilibrium  of a structure which occupies the open domain $\Omega\subset\mathbb{R}^d$ and which is submitted to  given body forces $f$,  to given traction forces $g$ on $\partial_f\Omega$  and to given displacements $u_0$ on the  complementary part  $\partial_u\Omega\neq\emptyset$. We  assume  the structure undergoes  small   perturbations  and  that  the  material   is  linear  elastic,
characterized by the Hooke's  tensor $\hooke$.  Let $u$ be  the unknown displacement field, $\varepsilon(u)$ the symmetric part  of the gradient, $\sigma$ the Cauchy stress tensor.

\begin{figure}[ht]\centering
  \includegraphics[width=.5\textwidth]{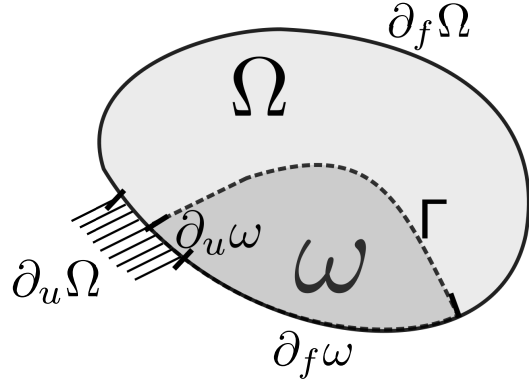}\caption{Domain $\Omega$, subdomain $\omega$ and boundaries}\label{fig:patate}
\end{figure}

Let $\omega  \subset \Omega$  be an  open subset of  $\Omega$, $\partial_f \omega=\partial \omega \cap \partial_f \Omega$, $\partial_u \omega=\partial \omega \cap \partial_u \Omega$ and $\Gamma=\partial\omega\setminus(\partial_u \omega\cup\partial_f \omega)$ (see Figure \ref{fig:patate}). We  introduce two
affine subspaces and one positive form:
\begin{itemize}
\item Subspace of kinematically admissible fields
\begin{equation}\label{eq:KA}
  \KA{\omega}=\left\{ u\in \left(\mathrm{H}^1(\omega)\right)^d,\ \trace (u) = u_0 \text{ on }\partial_u\omega \right\}
\end{equation}
where $\trace$ is the trace operator.
 \item Subspace of statically admissible fields
\begin{multline}\label{eq:SA}
  \SA{\omega}%, \hat{F},g,g}
  =\Bigg\lbrace   \tau\in  \left(\mathrm{L}^2(\omega)\right)^{d\times d}, %  \mathrm{H}(div, \omega),   
    \tau  \text{  symmetric}, \ %\hat{F}=\tau. n, \
    \forall  u^*\in  \KAoo{\omega},\ \\ \int_\omega
  \tau:\varepsilon(u^*)    d\omega    =    \int_\omega   f.u^*    d\omega +
  \int_{\partial_f\omega} g.u^* dS   \Bigg\rbrace
%  \int_{\partial_f\omega} {F}.u^* dS  \Bigg\rbrace
%  \int_{\partial\omega-\left(\partial_f\Omega\cup\partial_u\Omega\right)} \hat{F}.u^* dS  \Bigg\rbrace
\end{multline}
%where %$\KAoo{\omega}$ is the vector space associated to $\KA{\omega}$, 
\begin{equation*}
\text{where }  \KAoo{\omega}=\left\{ u\in \left(\mathrm{H}^1(\omega)\right)^d,\ \trace (u) = 0 \text{ on }\partial_u\omega\cup\Gamma \right\}
\end{equation*}
%
%$n$ the outward normal vector on $\partial \omega$ 
%and ${F}=g$ on $\partial_f \Omega$.
\item Measure of the non-verification of the constitutive equation \cite{ladeveze:75:erreur}
\begin{equation}\label{eq:ecr}
  \ecr{u,\sigma}{\omega} = \enernorm{\sigma-\hooke:\strain{u}}{\omega}
\end{equation}
where ${\enernorm{x}{\omega}}=\displaystyle \sqrt{\int_\omega \left( x: {\hooke}^{-1} :x \right)d\omega}$
\end{itemize}
The mechanical problem set on $\Omega$ can be formulated as:
\begin{center}
  Find    $\left(u_{ex},\sigma_{ex}\right)\in\KA{\Omega}\times\SA{\Omega}$    such    that
  $\ecr{u_{ex},\sigma_{ex}}{\Omega}=0$
\end{center}

%%%%%%%%%%%%%%%%%%%%%%%%%%%%%%%%%%%%%
%%%%%%%%%%%%%%%%%%%%%%%%%%%%%%%%%%%%% Discretization : Global problem
%%%%%%%%%%%%%%%%%%%%%%%%%%%%%%%%%%%%%
\subsection{Finite element approximation for the global problem}
Let  $\Omega_h$ be  a tessellation  of $\bar{\Omega}$  to which  we  associate a
finite  dimensional  subspace  $\KAh{\Omega}$  of $\KA{\Omega}$.  The  classical
finite element displacement approximation consists in searching
\begin{equation}
  \begin{aligned}
    u_h&\in\KAh{\Omega}\\
    \sigma_h&=\hooke:\varepsilon(u_h) \\ 
     \int_\Omega \sigma_h:\varepsilon(u_h^*)   d\Omega   &=   \int_\Omega  f.u_h^*   d\Omega  +  \int_{\partial_f\Omega} g.u_h^* dS, \qquad \forall u_h^*\in \KAhoo{\Omega}
  \end{aligned}
\end{equation}
%where
%\begin{multline}
%  \SAt{\omega}=\Bigg\lbrace\tau_h\in  \left(\mathrm{L}^2(\omega)\right)^{d\times
%    d}, \tau_h \text{  symmetric}, \   \forall u_h^*\in \KAhoo{\omega},\ \\ \int_\omega
%  \tau_h:\varepsilon(u_h^*)   d\omega   =   \int_\omega  f.u_h^*   d\omega   +
%  \int_{\partial_f\omega} g.u_h^* dS\Bigg\rbrace
%%  \int_{\partial_f\omega} {F_h}.u_h^* dS \Bigg\rbrace
%\end{multline}
%\begin{multline}
%  \SAt{\omega}=\Bigg\lbrace\tau_h\in  \left(\mathrm{L}^2(\omega)\right)^{n\times
%    n}, \tau_h \text{  symmetric}, \tau_h . n = \hat{F}_h, \   \forall u_h^*\in \KAho{\omega},\ \\ \int_\omega
%  \tau_h:\varepsilon(u_h^*)   d\omega   =   \int_\omega  g.u_h^*   d\omega   +
%  \int_{\partial\omega\cap\partial_f\Omega} g.u_h^* dS +
%  \int_{\partial\omega-\partial_f\Omega} \hat{F}_h.u^* dS \Bigg\rbrace
%\end{multline}
%Of course, except in trivial cases, we have $\sigma_h\notin\SA{\Omega}$.
%\begin{equation*}
%  \SA{\Omega}\subset \SAt{\Omega} \text{ and in most cases  }\sigma_h\notin\SA{\Omega}
%\end{equation*}
%$\sigma_h\notin\SA{\Omega}$.

After introducing the $d\times N_{dof}$ matrix $\shapev$ of  shape functions which form a basis of $\KAh{\Omega}$ and the vector of nodal unknowns  $\dep$ (of size $N_{dof}$, number of degrees of freedom) so that $u_h=\shapev \dep$, the  classical finite  element method leads to the well-known linear system: 
\begin{equation}\label{eq:globalFE}
  \stiff \dep = \force
\end{equation}
where $\stiff$ is  the (symmetric positive definite) stiffness  matrix of domain
$\Omega_h$ and $\force$ is the vector of generalized forces.

\subsection{A posteriori error estimator} %based on the error in constitutive relation}

The finite element approximation $(u_h,\sigma_h)$ satisfies $u_h\in\KA{\Omega}$ and $\ecr{u_h,\sigma_h}{\Omega}=0$ but $\sigma_h\notin\SA{\Omega}$. The error in constitutive relation consists in deducing from $(u_h,\sigma_h)$ an admissible displacement-stress pair
$(\admiss{u}_h,\admiss{\sigma}_h)\in\KA{\Omega}\times\SA{\Omega}$ in  order  to
measure the residual on the constitutive equation \eqref{eq:ecr} 
$\ecr{\admiss{u}_h,\admiss{\sigma}_h}{\Omega}\geqslant 0$. %(see \cite{ladeveze:04:mastCalculations} for more details). 
%Obviously, we have
%\begin{equation} \label{eq:crebasicscre2}
%  \ecr{\admiss{u}_h,\admiss{\sigma}_h}{\Omega}=0 \Longleftrightarrow (\admiss{u}_h,\admiss{\sigma}_h)=(u_{ex},\sigma_{ex})
%\end{equation}
%where $(u_{ex},\sigma_{ex})$ is the solution to the reference  problem.  
Using the well-known  Prager-Synge theorem it can be proved that
%always overestimates the true error:
\begin{equation*}
\strainnorm{\varepsilon(u_{ex})-\varepsilon(\admiss{u}_h)}{\structure}^2+\enernorm{\sigma_{ex}-\admiss{\sigma}_h}{\structure}^2 = \left(\ecr{\admiss{u}_h,\admiss{\sigma}_h}{\structure}\right)^2
\end{equation*}
Hence, the evaluation of the error in constitutive relation $\ecr{\admiss{u}_h,\admiss{\sigma}_h}{\structure}$ for any admissible pair ($\admiss{u}_h$, $\admiss{\sigma}_h$) provides a  guaranteed upper bound of the global error
\begin{equation}
  \label{eq:crebasicspragersynge}
\strainnorm{\varepsilon(u_{ex})-\varepsilon(\admiss{u}_h)}{\structure}\leqslant \ecr{\admiss{u}_h,\admiss{\sigma}_h}{\structure}
\end{equation}
%%%%%%%%%%%%%%%%%%%%%%%%%%%%%%%%%%%%%%%%%%%%%%%%%%%%%%%%%%%
%%%%%%%%%%%%%%%%%%%%%%%%%%%%%%%%%%%%%%%%%%%%%%%%%%%%%%%%%%%
%\subsection{Recovery admissible fields $(\admiss{u}_h,\admiss{\sigma}_h)$}
%\label{sec:crebasicssarecover}
%%%%%%%%%%%%%%%%%%%%%%%%%%%%%%%%%%%%%%%%%%%%%%%%%%%%%%%%%%%
%%%%%%%%%%%%%%%%%%%%%%%%%%%%%%%%%%%%%%%%%%%%%%%%%%%%%%%%%%%

$\KAh{\Omega}$ being a subspace of $\KA{\Omega}$, the construction of an admissible displacement field $\admiss{u}_h$ is straightforward since it can be taken equal to $u_h$. 
On the other hand, as $\sigma_h$ is not statically admissible, the construction of an admissible stress field $\admiss{\sigma}_h\in\SA{\Omega}$ is a crucial point which has already been widely studied in the literature. A first solution is to use a dual formulation of the reference problem \cite{beckers:1998} to compute $\admiss{\sigma}$ from scratch. Unfortunately building a subspace of $\SA{\Omega}$ is a complex task and most people prefer to post-process a statically admissible field from Field $\sigma_h$ obtained by a displacement formulation. 
Classical methods are the element equilibration techniques \cite{ladeveze:leguillon:81:erdc:sarecovery,ladeveze:rougeot:1997:erdc:sarecovery}, which have been improved by the use of the concept of partition of unity which lead to \cite{gallimard:2008:erdc:sarecovery,mointinhoDeAlmeida:maunder:2009:sarecoverypu,ladeveze:chamoin:florentin:2009:erdc:sarecovery} and the flux-free method \cite{pares:diez:huerta:2006:fluxfreesubdomain}. In most cases they involve the computation of efforts on ``star-patches'' which are the set of elements sharing one node, for each node of the mesh. Though rather simple these computations are in great number and thus expensive.

%\noindent
In the following, we note by $\mathcal{F}_h$ the algorithm which has been chosen to build an admissible stress field $\hat{\sigma}_h$. Whatever the choice, the algorithm takes as input not only the finite element stress field $\sigma_h$ but also the continuous representation of the imposed forces $(f, g)$.
\begin{equation*}
\hat{\sigma}_h=\mathcal{F}_h(\sigma_h, f, g) \in \SA{\Omega} %\, ; \, \sigma_h \in \SAt{\Omega}.
\end{equation*}
The algorithm we have used for our applications is the one proposed in \cite{ladeveze:rougeot:1997:erdc:sarecovery} using a three degrees higher polynomial basis when solving the local problems on elements \cite{babuska:94:vpeena}.

% Since  $\KAh{\Omega}$ is  a  subspace of  $\KA{\Omega}$,  $\admiss{u}_h$ can  be
% chosen equal to $u_h$.  However  $\sigma_h$ is not statically admissible so that
% $\admiss{\sigma}_h\in\SA{\Omega}$ has to  be deduced from $\sigma_h$.  Different
% approaches   to   recover    $\admiss{\sigma}_h$   exist   in   the   literature
% \cite{ladeveze:leguillon:81:erdc:sarecovery,ladeveze:rougeot:1997:erdc:sarecovery},
% some     of    which     based     on    partition     of    unity     principle
% \cite{gallimard:2008:erdc:sarecovery,
%   mointinhoDeAlmeida:maunder:2009:sarecoverypu,
%   ladeveze:chamoin:florentin:2009:erdc:sarecovery}, and we choose to use the one
% proposed by \cite{ladeveze:leguillon:81:erdc:sarecovery}.

%Note that the used of patches can not be employed on the boundary nodes without assuming communication between subdomains. Though these exchanges would remain limited, we propose an alternate strategy to achieve full parallelism without impairing the properties of the error in constitutive relation.
%In the  following, we  show how these  two difficulties  are dealt with  for two amongst  the most  famous  domain decomposition  strategies,  namely the  primal approach (BDD) and the dual approach  (FETI); we make use of results specifically described in \cite{RIXEN:1998:SUPERL,GOSSELET:2003:IEI,GOSSELET.2007.1}.
%%%%%%%%%%%%%%%%%%%%%%%%%%%%%%%%%%%%%%
%%%%%%%%%%%%%%%%%%%%%%%%%%%%%%%%%%%%%%
%%%%%%%%%%%%%%%%%%%%%%%%%%%%%%%%%%%%%%
\subsection{Substructured formulation}
%%%%%%%%%%%%%%%%%%%%%%%%%%%%%%%%%%%%%%
%%%%%%%%%%%%%%%%%%%%%%%%%%%%%%%%%%%%%%
%%%%%%%%%%%%%%%%%%%%%%%%%%%%%%%%%%%%%%
%We consider a non overlapping domain decomposition  of $\Omega=\displaystyle \cup_s\Omega\s$ ($\Omega\s\cap\Omega^{(s')}=\emptyset$ for $s\neq s'$). %and $\bar{\Omega}=\cup_s  \bar{\Omega}\s$. 
Let us consider a decomposition  of domain $\Omega$ in open subsets $(\Omega\s)_{1\leqslant s\leqslant N_{sd}}$ ($N_{sd}$ is the number of subdomains) so that   $\Omega\s\cap\Omega^{(s')}=\emptyset$    for    $s\neq   s'$    and $\bar{\Omega}=\cup_s  \bar{\Omega}\s$.
Let $u\g=(u\s)_s$,  we define  the global assembling operator $\assemg$:
\begin{equation}
\begin{aligned}
  u= \assemg (u\g) \Leftrightarrow u_{|\Omega\s}=u\s
\end{aligned}
\end{equation}
In order to reformulate the mechanical problem on the substructured configuration,
we need  to specify the  conditions that  should be satisfied  at the
boundary between subdomains $\Gamma^{(ss')}=\partial\Omega\s\cap\partial\Omega^{(s')}$.
%(which does not enter in the description of previous admissible subspaces). 
We have the fundamental properties:
\begin{equation}\label{eq:KAss}
  \assemg (u\g) \in \KA{\Omega} \Leftrightarrow \left\lbrace \begin{array}{l} u\s\in \KA{\Omega\s},\ \forall s \\ \trace(u\s)=\trace(u^{(s')}) \text{ on } \Gamma^{(ss')},\ \forall(s,s')\end{array}\right.
\end{equation}
\begin{equation}\label{eq:SAss}
  \assemg (\sigma\g) \in \SA{\Omega} \Leftrightarrow \left\lbrace \begin{array}{l} \sigma\s\in \SA{\Omega\s},\ \forall s \\ 
%{F}^{(s)} + {F}^{(s')}=0 \text{ on } \Gamma^{ss'},\ \forall(s,s')\end{array}\right.
\sigma\s.n\s+\sigma^{(s')}.n^{(s')}=0 \text{ on }\Gamma^{(ss')},\ \forall(s,s')\end{array}\right.
%\int_{\partial\Omega\s\cap\partial\Omega^{(s')}} \hat{F}\s.u^* dS+\int_{\partial\Omega\s\cap\partial\Omega^{(s')}} \hat{F}\s.u^* dS =0 \ ,\ \forall(s,s')\end{array}\right.
\end{equation}
%where $\Gamma^{ss'}=\partial\Omega\s\cap\partial\Omega^{(s')}$.\\
% au sens des normes qui vont bien ?
%where $n\s$ stand  for the outer normal to $\Omega\s$. 
In  other words, in order to  be admissible  on the  whole domain  $\Omega$, not  only fields  need  to be
admissible in a local sense (independently on  each $\Omega\s$), but they also need to satisfy interface  conditions,  namely displacements continuity and tractions balance (action-reaction principle).
%$$
%\int_{\partial\Omega\s\cap\partial\Omega^{(s')}} (\hat{F}\s + \hat{F}\s).u^* dS =0 \ \forall u^* \in \KAo{\Omega}
% \ , \  \forall(s,s')
%$$

%%%%%%%%%%%%%%%%%%%%%%%%%%%%%%%%%%%%%
%%%%%%%%%%%%%%%%%%%%%%%%%%%%%%%%%%%%% Discretization : substructured context
%%%%%%%%%%%%%%%%%%%%%%%%%%%%%%%%%%%%%
\subsection{Finite element approximation for the substructured problem}
\begin{figure}[H]\centering
\psfrag{a}{$1^{(1)}$}\psfrag{b}{$2^{(1)}$}\psfrag{c}{$3^{(1)}$}\psfrag{d}{$4^{(1)}$}\psfrag{e}{$5^{(1)}$}
\psfrag{f}{$1^{(2)}$}\psfrag{g}{$2^{(2)}$}\psfrag{h}{$3^{(2)}$}\psfrag{i}{$4^{(2)}$}\psfrag{j}{$5^{(2)}$}
\psfrag{k}{$1^{(3)}$}\psfrag{l}{$2^{(3)}$}\psfrag{m}{$3^{(3)}$}\psfrag{n}{$4^{(3)}$}
\psfrag{ct}{$1_b^{(1)}$}\psfrag{dt}{$2_b^{(1)}$}\psfrag{et}{$3_b^{(1)}$}
\psfrag{ht}{$1_b^{(2)}$}\psfrag{it}{$2_b^{(2)}$}\psfrag{jt}{$3_b^{(2)}$}
\psfrag{kt}{$1_b^{(3)}$}\psfrag{lt}{$2_b^{(3)}$}\psfrag{mt}{$3_b^{(3)}$}
\psfrag{o}{$1_\Gamma$}\psfrag{p}{$2_\Gamma$}\psfrag{q}{$3_\Gamma$}\psfrag{r}{$4_\Gamma$}
\psfrag{s}{$\underline{1}_\Gamma$}\psfrag{t}{$\underline{2}_\Gamma$}\psfrag{u}{$\underline{3}_\Gamma$}\psfrag{v}{$\underline{4}_\Gamma$}\psfrag{w}{$\underline{5}_\Gamma$}\psfrag{x}{$\underline{6}_\Gamma$}
\subfigure[Subdomains]{\includegraphics[width=5.cm]{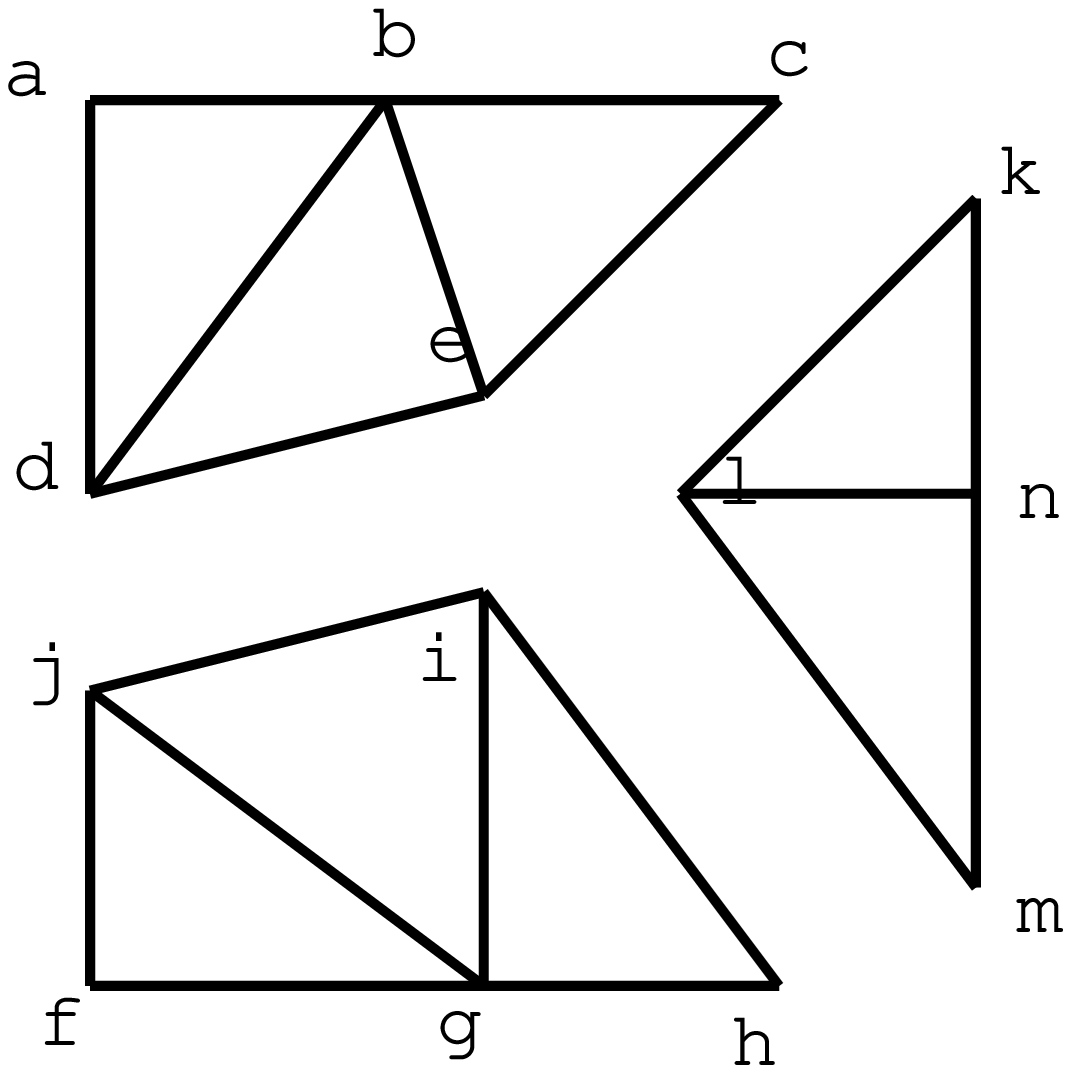}
}\qquad
\subfigure[Local interface]{\includegraphics[width=5.cm]{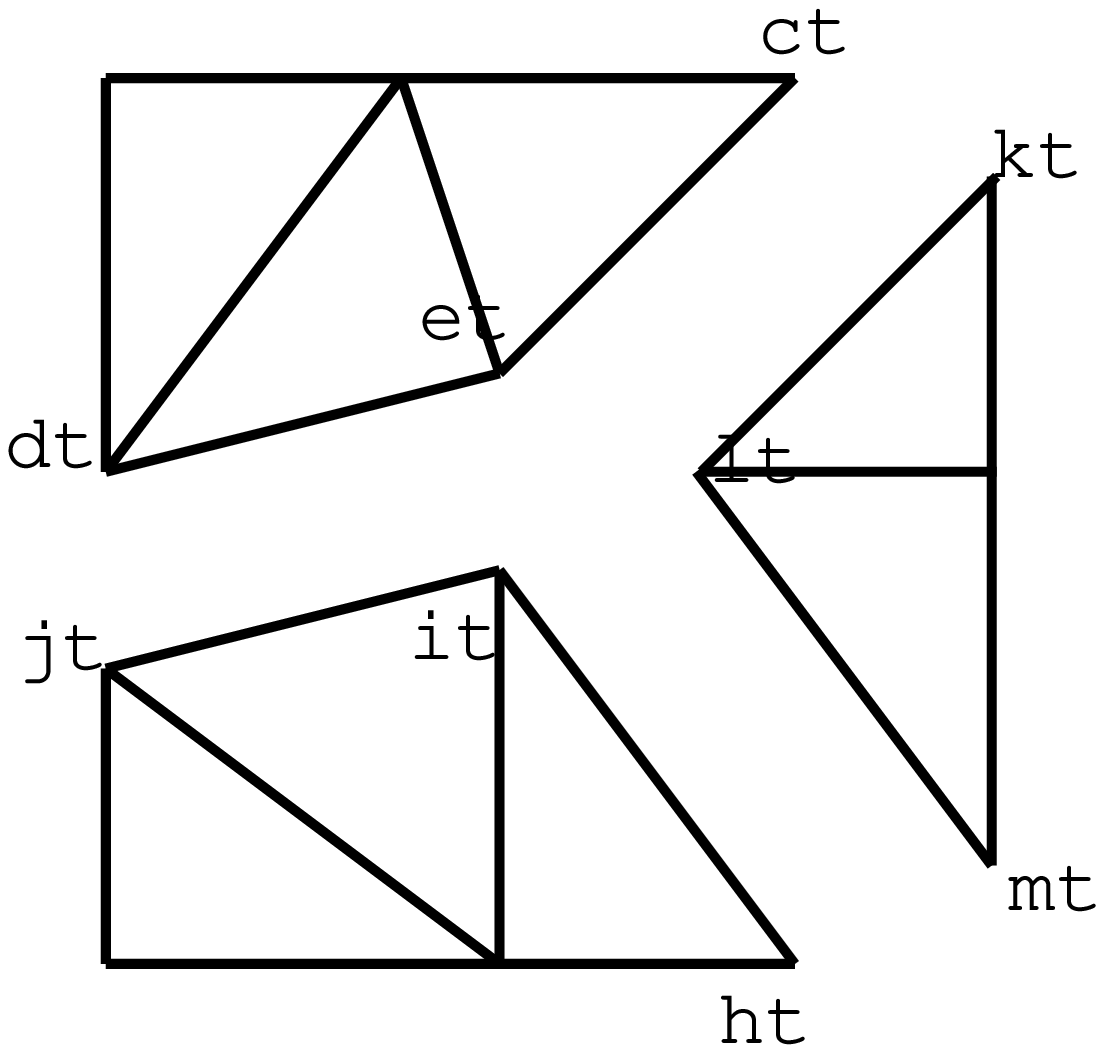} }\\
\subfigure[Primal interface]{\includegraphics[width=5.cm]{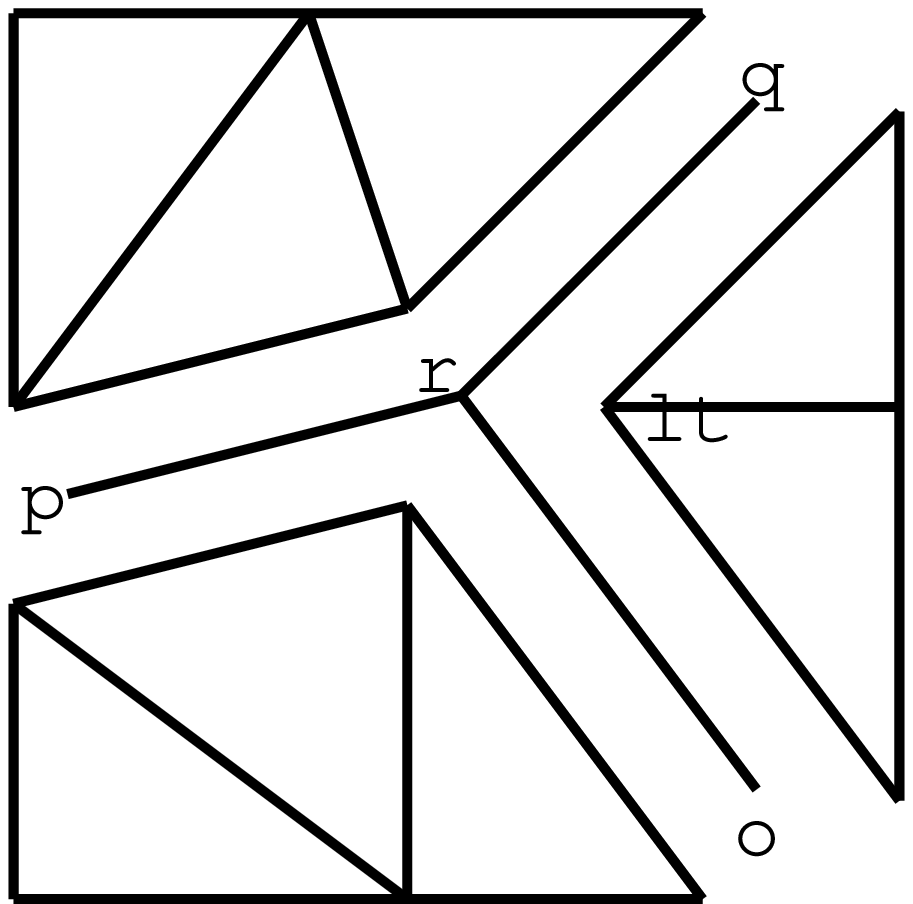}
}\qquad \subfigure[Dual interface]{\includegraphics[width=5.cm]{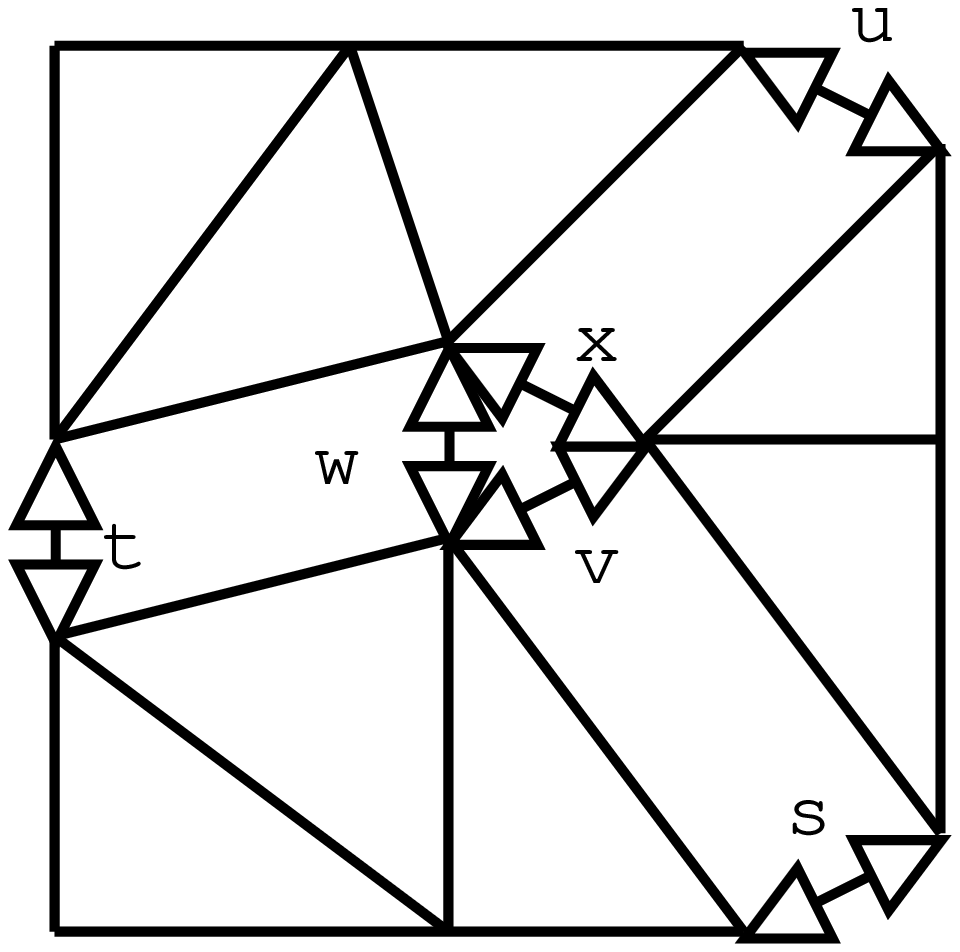}
} \scriptsize\begin{equation*}
\begin{array}{rclrclrcl}
\traceh^{(1)}&=&\begin{pmatrix}0&0&1&0&0\\0&0&0&1&0\\0&0&0&0&1\end{pmatrix},&\traceh^{(2)}&=&\begin{pmatrix}0&0&1&0&0\\0&0&0&1&0\\0&0&0&0&1\end{pmatrix},&\traceh^{(3)}&=&\begin{pmatrix}1&0&0&0\\0&1&0&0\\0&0&1&0\end{pmatrix}\\
\passem^{(1)}&=&\begin{pmatrix}0&0&0\\0&1&0\\1&0&0\\0&0&1\end{pmatrix},&\passem^{(2)}&=&\begin{pmatrix}1&0&0\\0&0&1\\0&0&0\\0&1&0\end{pmatrix},&\passem^{(3)}&=&\begin{pmatrix}0&0&1\\0&0&0\\1&0&0\\0&1&0\end{pmatrix}\\
\dassem^{(1)}&=&\begin{pmatrix}0&0&0\\0&1&0\\1&0&0\\0&0&0\\0&0&1\\0&0&1\end{pmatrix},&\dassem^{(2)}&=&\begin{pmatrix}1&0&0\\0&0&-1\\0&0&0\\0&1&0\\0&-1&0\\0&0&0\end{pmatrix},&\dassem^{(3)}&=&\begin{pmatrix}0&0&-1\\0&0&0\\-1&0&0\\0&-1&0\\0&0&0\\0&-1&0\end{pmatrix}
\end{array}\end{equation*}\normalsize
\caption{Local numberings, interface numberings, trace and assembly operators}\label{fig:omegef:2}
\end{figure}

We assume that the tessellation of  $\bar{\Omega}$ and the substructuring are conforming so that (i) each  element only belongs to one  subdomain and (ii) nodes are  matching on the interfaces.   Each  degree of  freedom  is  either  located inside  a  subdomain
(subscript $i$)  or on its boundary $\Gamma\s=\cup_{s'}\Gamma^{(ss')}$ (subscript  $b$) where it is  shared with at
least one neighboring subdomain. Let $\lam_b\s$ be the vector of unknown efforts
imposed on the interface of  subdomain $\Omega\s_h$ by its neighbors. The finite
element   problem   \eqref{eq:globalFE}   can   be  written   highlighting   the
contributions of subdomains:
\begin{equation}\label{eq:inter_bool}
\forall s,\ \stiff\s\dep\s= \force\s + {\traceh\s}^T\lam_b\s \text{ with }
\left\{\begin{array}{l}
\sum\limits_s \passem\s \lam_b\s = \mathbf{0}\\
\sum\limits_s \dassem\s \dep_b\s =  \mathbf{0}
\end{array}\right.
\end{equation}
where  $\traceh\s$ is the  discrete trace  operator ($\dep\s_b=\traceh\s\dep\s$)
and  where  $\passem\s$  and   $\dassem\s$  are  assembling  operators  so  that
$\passem\s$  enables  to  formulate  the mechanical  equilibrium  of  interfaces
\eqref{eq:SAss}  and   $\dassem\s$  enables  to  formulate   the  continuity  of
displacements  \eqref{eq:KAss}   (in  the  case  of  two   subdomains,  we  have
$\sum_s\passem\s\lam_b=\lam_b^{(1)}+\lam_b^{(2)}=\mathbf{0}$                  and
$\sum_s\dassem\s\dep_b\s=\dep_b^{(1)}-\dep_b^{(2)}=\mathbf{0}$, see Fig.~\ref{fig:omegef:2} for less trivial example and \cite{GOSSELET.2007.1} for more an extensive description of all operators). One fundamental
property of assembling operators is their orthogonality:% \cite{GOSSELET.2007.1}:
\begin{equation}\label{ortho_assem}
\sum_s \dassem\s{\passem\s}^T=\mathbf{0}
\end{equation}
Note that the equilibrium of subdomain $\Omega\s$ also writes: 
\begin{equation}\label{eq:equi_sd}
 \begin{pmatrix} \stiff_{ii}\s & \stiff_{ib}\s \\ \stiff_{bi}\s & \stiff_{bb}\s \end{pmatrix}
\begin{pmatrix} \dep_{i}\s \\ \dep_{b}\s \end{pmatrix} = \begin{pmatrix} \force_{i}\s \\ \force_{b}\s \end{pmatrix} +\begin{pmatrix} \mathbf{0}_{i}\s \\ \lam_{b}\s \end{pmatrix}
\end{equation}
or in an equivalent condensed form:
\begin{equation}\label{eq:equi_sd_cond}
\schur\s \dep_{b}\s   = \rhs_p\s + \lam_{b}\s
\end{equation}
with
\begin{equation}\label{eq:oper_sd_cond}
\begin{aligned}
\schur\s & =\stiff\s_{bb}-\stiff\s_{bi}{\stiff\s_{ii}}^{-1}\stiff\s_{ib}\\
\rhs\s & =\force_b\s -\stiff\s_{bi}{\stiff\s_{ii}}^{-1}\force_i\s\\
\end{aligned}
\end{equation}
where  $\schur\s$ is the Schur complement and $\rhs\s$ is the condensed right-hand side.

%%%%%%%%%%%%%%%%%%%%%%%%%%%%%%%%%%%%%%%%%%%%%%%%%%%%%%%%%%%%%%%%%%%%%
%%%%%%%%%%%%%%%%%%%%%%%%%%%%%%%%%%%%%%%%%%%%%%%%%%%%%%%%%%%%%%%%%%%%%
%%%%%%%%%%%%%%%%%%%%%%%%%%%%%%%%%%%%%%%%%%%%%%%%%%%%%%%%%%%%%%%%%%%%%
%%%%%%%%%%%%%%%%%%%%%%%%%%%%%%%%%%%%%%%%%%%%%%%%%%%%%%%%%%%%%%%%%%%%%
\section{A posteriori error estimator in substructured context}\label{sec:errordd} %based on the error in constitutive relation}

%The key point for an efficient evaluation of the error in constitutive relation $\ecr{\admiss{u}_h,\admiss{\sigma}_h}{\structure}$ in a substructured context (without overlap), is to define admissible pairs ($\admiss{u}\s_h$,$\admiss{\sigma}\s_h$) $\in \KA{\Omega\s} \times\SA{\Omega\s}$ on each subdomain so that the associated assembled pair $(\assemg (\admiss{u}_h\g),\assemg (\admiss{\sigma}_h\g))=(\admiss{u}_h,\admiss{\sigma}_h)$ is an admissible pair for the reference problem (that is an element of $\KA{\Omega}\times\SA{\Omega}$). Due to the absence of overlap, the additive structure of the associated error in constitutive relation
The key point for the efficient evaluation of the error in constitutive relation in a substructured context (without overlapping) is to define admissible pairs ($\admiss{u}\s_h$,$\admiss{\sigma}\s_h$) $\in \KA{\Omega\s} \times\SA{\Omega\s}$ on each subdomain so that the associated assembled pair is admissible for the reference problem $(\assemg (\admiss{u}_h\g),\assemg (\admiss{\sigma}_h\g))\in\KA{\Omega}\times\SA{\Omega}$. Due to the absence of overlap, the additive structure of the associated error in constitutive relation leads to a fully parallel evaluation of the a posteriori error estimator:
\begin{equation*}
%\ecr{\admiss{u}_h,\admiss{\sigma}_h}{\structure} =\ecr{\assemg (\admiss{u}_h\g),\assemg (\admiss{\sigma}_h\g)}{\structure} = \sum_s \ecr{\admiss{u}\s_h,\admiss{\sigma}\s_h}{\Omega\s}
\left(\ecr{\assemg (\admiss{u}_h\g),\assemg (\admiss{\sigma}_h\g)}{\structure}\right)^2 = \sum_s \left(\ecr{\admiss{u}\s_h,\admiss{\sigma}\s_h}{\Omega\s}\right)^2
\end{equation*}
%leads then to a fully parallel evaluation of the a posteriori error estimator. 

The application of a  classical  recovery   strategy  to compute admissible fields raises  two  difficulties  in  a substructured context.
%\begin{itemize}
%\item  
First, the  star-patches can  not  be employed  on  the  boundary nodes  without
  assuming communication between subdomains. Though these exchanges would remain
  limited, we propose an alternate  strategy to achieve full parallelism without
  impairing the properties of the error in constitutive relation.
%\item  
Second, in  order  to   solve  the  substructured  problem  \eqref{eq:inter_bool}
 parallel  strategies consist in  using iterative solvers which are  based on  the loosening of at least  one of the interface  conditions which is only  verified (up to  a certain  precision) once  the solver  converged. Thus recovering strategies need to be adapted so that the local fields $(\admiss{u}_h\s,\admiss{\sigma}_h\s)$ satisfy the interface conditions.
%\end{itemize}

The aim of this section is to prove that the determination of the admissible pair $(\assemg(\admiss{u}_h\g),\assemg (\admiss{\sigma}_h\g))$ can be brought back to the determination of nodal interface fields $(\hdep_b\s,\hlam_b\s)_s$ which satisfy specific interface conditions. The construction of these nodal fields depends on the chosen domain decomposition strategy and is discussed in the following sections.

\subsection{Kinematically admissible fields}
In order to ensure interface Condition \eqref{eq:KAss} when building $\admiss{u}\s_h \in \KA{\Omega\s}$ so that $\assemg (\admiss{u}_h\g) \in \KA{\Omega}$, we introduce continuous interface displacement fields $\hat{u}_{bh}\s$ from which we shall deduce internal displacement fields:
\begin{equation*}
\begin{aligned}
\hat{u}\s_{bh}&= \hat{u}^{(s')}_{bh} \ , \ \forall (s,s')\\
\hat{u}\s_{h|\Gamma^{(ss')}} &= \hat{u}\s_{bh}, \ \forall s
\end{aligned}
\end{equation*}
Since discretizations are matching on the interface, the first condition can directly be imposed on finite element nodal quantities:
\begin{equation*}{\hdep_b\s} = {\hdep^{(s')}_b}\ , \ \forall (s,s')\end{equation*}
In order to deduce the internal fields, one finite element problem is solved independently on each subdomain with imposed Dirichlet conditions on the interface:
\begin{equation*}
\begin{aligned}
  {\hdep_i\s}&={\stiff\s_{ii}}^{-1}\left(\force_i\s-\stiff\s_{ib}{\passem\s}^T\hdep_b\s\right)\\
  \hat{u}_h\s&=\shapev\s \hdep\s= \begin{pmatrix}{\shapev}_{_i}\s& {\shapev}_{_b}\s\end{pmatrix}\begin{pmatrix} \hdep_i\s \\\hdep_b\s  \end{pmatrix}\\
\admiss{u}&=\assemg (\admiss{u}_h\g) \in \KA{\Omega}
\end{aligned}
\end{equation*}

%The construction of an admissible displacement field $\admiss{u}\s_h \in \KA{\Omega\s}$ so that $\assemg (\admiss{u}_h\g) \in \KA{\Omega}$ only requires to impose
%on $\Gamma^{ss'}=\partial \Omega\s\cap\partial \Omega^{(s')}$,
%\begin{equation*}u_{bh}={u\s_h}_{|\Gamma^{ss'}}= {u^{(s')}_h}_{|\Gamma^{ss'}} \ , \ \forall (s,s')\end{equation*}
%The choice and the construction of $u_{bh}$ depend on the domain decomposition method.
%which can be done very easily (see section).\\
\subsection{Statically admissible fields}
In order to ensure interface Condition \eqref{eq:SAss} when building $\admiss{\sigma}\s_h \in \SA{\Omega\s}$ so that $\assemg (\admiss{\sigma}_h\g) \in \SA{\Omega}$, we introduce for each subdomain the continuous balanced interface traction fields $\adF_{bh}\s$ defined on $\Gamma\s$ which satisfy:
\begin{equation}\label{eq:Fadm}
\begin{aligned}
\admiss{\sigma}_h\s.n\s &= \adF_{bh}\s \text{ on }\Gamma\s\\
\adF_{bh}\s + \adF_{bh}^{(s')} &=0 \text{ on }\Gamma^{(ss')}\\ 
\int_{\Omega\s} f .\rho d\Omega + \int_{\partial_f\Omega\s} g\s .\rho dS + \int_{\Gamma\s} \adF_{bh}\s .\rho dS & =0\quad \forall\rho \in \mathrm{RKA}^0(\Omega\s)
\end{aligned}
\end{equation}
where $\mathrm{RKA}^0(\Omega\s)$ is the set of rigid body motions which are compatible with Dirichlet conditions imposed on $\partial_u\Omega\s$:% and Neumann conditions imposed on  $\Gamma\s\cup\partial_f\Omega\s$:
\begin{equation*}
  \mathrm{RKA}^0(\Omega\s)=\left\{ \rho\in\mathrm{H}^1(\Omega\s),\ \rho=0 \text{ on }\partial_u\Omega\s,\ \varepsilon(\rho)=0, \right\}
\end{equation*}
The last condition of \eqref{eq:Fadm} is the translation of Fredholm's alternative in order to ensure the well-posedness of the static problem on domain $\Omega\s$.
To build these traction fields in a simple way, we associate them with the finite element nodal reaction field ${\hlam}\s_b$:
\begin{equation}\label{eq:salambda1}
    \int_{\Gamma^{(ss')}} \adF_{bh}\s\cdot{\shapef_j\s}_{|\Gamma^{(ss')}}dS = {\hlam}_{b,j}\s
\end{equation}
where $j$ denote a node of the interface, $\shapef_j\s$  its associated shape  function and ${\hlam}\s_{b,j}$ the corresponding nodal component of ${\hlam}\s_b$. This equation then imposes that the discrete field ${\hlam}_{b}\s$ and the continuous field $\adF_{bh}\s$ develop the same virtual work in any finite element displacement field. The conditions on $\adF_{bh}\s$ have these discrete counterparts on ${\hlam}\s_b$:
\begin{equation}%\label{eq:salambda2}
  \begin{aligned}
    \sum\limits_s \passem\s \hlam_b\s &= \mathbf{0}\\
    {\kernel\s}^T\left({\traceh\s}^T{\hlam}_{b}\s+\force\s\right)&=\mathbf{0}
  \end{aligned}
\end{equation}
where  $\kernel\s$ is a basis of $\ker(\stiff\s)$.
 As said earlier, the first equation corresponds to the equilibrium between subdomains. The second equation corresponds to the balance of the subdomain with respect to virtual rigid body motions (since this kind of displacement field is exactly represented in the finite element approximation, the discrete condition is equivalent to the continuous one).

%Thus starting from any $(\lam_b\s)_s$ which satisfies \eqref{eq:salambda2}, any continuous field satisfying \eqref{eq:salambda1} will enable to build independently on each $\Omega\s$ an admissible field $\admiss{\sigma}_h\s$  so that $\assemg (\admiss{\sigma}_h\g) \in \SA{\Omega}$.
%we deduce them from balanced nodal interface forces ${\lam}\s_b$ (with $\sum\limits_s \passem\s \lam_b\s = \mathbf{0}$). The constitency of virtual work imposes
%The construction of an admissible stress field $\admiss{\sigma}\s_h \in \SA{\Omega\s}$ so that $\assemg (\admiss{\sigma}_h\g) \in \SA{\Omega}$ requires to impose:
%\begin{equation*}\int_{\Gamma^{ss'}} (\admiss{\sigma}_h\s.n\s +\admiss{\sigma}_h^{(s')}.n^{(s')}).u^* dS = 0 \ ,\ \forall u^* \in \KAoo{\Omega}\end{equation*}
%which will be satisfy if we impose on each finite element segment of on each $\Gamma_h^{ss'}}$
%$${F}_h^{(s)} + {F}_h^{(s')}=0  \, on \Gamma_h^{ss'} \ \forall(s,s').$$
%To achieve this last condition, a simple way consists in defining continuous balanced boundary conditions $F_{h}\s$ on each interface.
%We choose to built $F_{h}\s$ from a balanced nodal interfaces forces ${\lam}\s_b$ ($\sum\limits_s \passem\s \lam_b\s = \mathbf{0}$) according to 
As a first approach, we define $\adF_{bh}\s$ as:
\begin{equation}
\label{eq:salambda3}
\adF_{bh}\s = \shapev{}_{|\Gamma^{(s)}}\s {\hF}_{b}\s
\end{equation}
where ${\hF}_{b}\s$ is the vector of nodal values of $\adF_{bh}\s$ and $\shapev{}_{|\Gamma^{(s)}}\s$ refers to the vector of the trace on $\Gamma^{(s)}$ of finite element shape functions.
Vector ${\hF}_{b}\s$ is then obtained by the inversion of the (small) ``mass'' matrix of the interface of each subdomain. 
%Besides, if nodal interface forces $\lam_{b}\s$ come from a well-posed discrete Neumann problem on subdomain $\Omega\s$, the problem on subdomain $\Omega\s$ with Neumann boundary condition ${F}_{h}\s$ on $\partial \Omega\s$ is also well-posed.
In the following, we denote by $\mathcal{G}_h$ the previous procedure which associates a continuous balanced interface force $\adF_{bh}\s$ to a balanced nodal interfaces forces ${\hlam}\s_b$:
\begin{equation}
\label{eq:Fh}
\adF_{bh}\s=\mathcal{G}_h({\hlam}\s_b)
\end{equation}

The traction field $\adF_{bh}\s$ allows to satisfy the interface conditions associated to the static admissibility. The next step is to build internal finite element stress fields which match the associated nodal boundary field ${\hlam}\s_b$. This is done by solving one finite element problem on each subdomain with imposed Neumann conditions on the interface.
\begin{equation}\label{eq:bdd_pi}
\begin{aligned}
\tilde{\dep}\s &={\stiff\s}^+\left(\force\s+{\traceh\s}^T{\hlam}\s_b\right)\\
%\tilde{\sigma}_h\s &= \hooke:\varepsilon(\shapev\s\tilde{\dep}\s) \in \SAt{\Omega\s, {F}_h\s}\\
\admiss{\sigma}_h\s& =\mathcal{F}_h\left(\hooke:\varepsilon(\shapev\s\tilde{\dep}\s), f\s, \left\{g\s,\mathcal{G}_h({\hlam}\s_b)\right\}\right)\\
\admiss{\sigma}_h&=\assemg (\admiss{\sigma}_h\g) \in \SA{\Omega}
\end{aligned}
\end{equation}
The use  of the pseudo-inverse  ${\stiff\s}^+$ is due  to the potential lack of Dirichlet boundary conditions on the substructure. Displacement field $\tilde{\dep}\s$ is defined up to a rigid body motion which needs not to be determined since only the symmetric gradient of the associated displacement field is required.

It has to be noted that the fully parallel procedure $\mathcal{G}_h$ proposed above leads to a different admissible traction field as would have been obtained using standard patch-technique \cite{ladeveze:rougeot:1997:erdc:sarecovery} (referred in the sequel as the sequential approach). Thus the use of $\mathcal{G}_h$ implies that the parallel error estimation is different from the standard sequential one even when discrete interface conditions are satisfied. For now there are no theoretical results on the quality of the resulting fields, examples (as given in Section~\ref{sec:numeric}) show that sequential estimator and parallel estimator (when interface conditions have sufficiently converged, which happens very quickly) can not be distinguished.

%The construction of admissible stress field $\admiss{\sigma}\s_h \in \SA{\Omega\s}$ such that $\assemg (\admiss{\sigma}_h\g) \in \SA{\Omega}$ is finally reduced to the construction of balanced nodal interface forces $\lam_{b}\s$ associated to well-posed discrete Neumann problem. Anyhow, the construction and the choice of $\lam_{b}\s$ depends on the domain decomposition method.

%%%%%%%%%%%%%%%%%%%%%%%%%%%%%%%%%%%%%%%%%%%%%%%%%%%%%%%%%%%%%%%%%%%%%
%%%%%%%%%%%%%%%%%%%%%%%%%%%%%%%%%%%%%%%%%%%%%%%%%%%%%%%%%%%%%%%%%%%%%
%%%%%%%%%%%%%%%%%%%%%%%%%%%%%%%%%%%%%%%%%%%%%%%%%%%%%%%%%%%%%%%%%%%%%
%%%%%%%%%%%%%%%%%%%%%%%%%%%%%%%%%%%%%%%%%%%%%%%%%%%%%%%%%%%%%%%%%%%%%

\section{Recovery of admissible fields in BDD}
\label{sec:bddrecovery}
In the  Balancing domain decomposition \cite{MANDEL:1993:BAL,LETALLEC:1994:DDM},
a  unique  interface  displacement   unknown  $\dep_b$  is  introduced  so  that
continuity is always insured:
\begin{equation}
  \dep_b\s={\passem\s}^T \dep_b \Longrightarrow \sum_s\dassem\s\dep_b\s = \mathbf{0}
\end{equation}
Other quantities can be deduced from $\dep_b$ and equations (\ref{eq:equi_sd},\ref{eq:equi_sd_cond}):
\begin{equation}\label{eq:bdd1}
\begin{aligned}
\dep_i\s &= {\stiff\s_{ii}}^{-1}\left(\force_i\s-\stiff\s_{ib}{\passem\s}^T\dep_b\right)\\
\lam_b\s &= \schur\s \dep_b - \rhs\s
\end{aligned}
\end{equation}
%\begin{equation}\label{eq:bdd1}
%\begin{aligned}
%\dep_i\s &= {\stiff\s_{ii}}^{-1}\left(\force_i\s-\stiff\s_{ib}{\passem\s}^T\dep_b\right)\\
%\lam_b\s &= \left(\stiff\s_{bb}-\stiff\s_{bi}{\stiff\s_{ii}}^{-1}\stiff\s_{ib}\right){\passem\s}^T\dep_b - \left(\force_b\s %-\stiff\s_{bi}{\stiff\s_{ii}}^{-1}\force_i\s\right). 
%\end{aligned}
%\end{equation}
%where         we          recognize         the         Schur         complement
%$\schur\s=\left(\stiff\s_{bb}-\stiff\s_{bi}{\stiff\s_{ii}}^{-1}\stiff\s_{ib}\right)$
%and      the     condensed     right-hand      Side     $\rhs\s=\left(\force_b\s
%  -\stiff\s_{bi}{\stiff\s_{ii}}^{-1}\force_i\s\right)$.  
%The resolution consists in iteratively trying to verify the last condition:
%\begin{equation}\label{eq:primal_sys}
%\mathbf{0}=\sum_s \passem\s\lam_b\s
%= \left(\sum_s \passem\s\schur\s{\passem\s}^T\right)\dep_b - \left(\sum_s \passem\s\rhs_p\s\right)
%\end{equation}
The BDD solver consists in iteratively finding the interface displacement $\dep_b$ which insure global equilibrium
($\sum_s\passem^{(s)}\lam^{(s)}_b=\mathbf{0}$),
\begin{equation}\label{eq:primal_sys}
\mathbf{0}=\sum_s \passem\s\lam_b\s
= \left(\sum_s \passem\s\schur\s{\passem\s}^T \right)\dep_b - \left(\sum_s \passem\s\rhs\s\right)
%= \schur \dep_b - \rhs
\end{equation}
%with  
%$$
%\schur = \sum_s \passem\s\schur\s{\passem\s}^T \ ; \ \rhs = \sum_s \passem\s\rhs\s
%$$

\subsection{Recovery of KA fields}
In the BDD solver, kinematic interface conditions are satisfied anytime and using
$\hdep_b\s=\dep_b$ enables to build  $\admiss{u}\s_h$ so that $\admiss{u}_h=\assemg\left(\admiss{u}_h\g\right)\in\KA{\Omega}$. Note that all associated 
computations are realized  during the standard resolution process so that no extra operation is required.

\subsection{Recovery of SA fields}
For a given interface displacement $\dep_b$, we note:
\begin{equation*}
\llceil \lam_b \rrceil=\sum_s\passem^{(s)}\lam^{(s)}_b=\sum_s\passem^{(s)} \left( \schur\s \dep_b - \rhs\s \right)%=\schur \dep_b - \rhs
\end{equation*}
Obviously $\llceil \lam_b \rrceil$ is zero if and only if $\dep_b$ is the solution to \eqref{eq:primal_sys}. We then define:
\begin{equation}\label{eq:lam_feti}
{\hlam}\s_b  = \lam\s_b- \left.\psassem\s\right.^T\llceil \lam_b \rrceil
\end{equation}
where $(\psassem\s)_s$ are scaled  assembling  operators so that  $\sum_s\passem\s\left.\psassem\s\right.^T= \mathbf{I}$. The multiplicity scaling is a typical example of such operator $\psassem\s$:
\begin{equation*}
\left.\psassem\s\right.^T = \left.\passem\s\right.^T\left(\sum_j\passem^{(j)}\left.\passem^{(j)}\right.^T\right)^{-1}
\end{equation*}
which, in  the  case  of  two subdomains, gives $\left.\psassem\s\right.^T \llceil \lam_b \rrceil=\displaystyle \frac{1}{2}\llceil \lam_b \rrceil$. In the case of heterogeneous structures, other scaled assembly operators which take the heterogeneity into account are used \cite{RIXEN:1998:SUPERL,KLAWONN:2001:FNN,GOSSELET:2002:DDM}.

It is clear that by definition, $\hlam\s_b$ is a balanced nodal reaction field:
%\begin{equation}\label{eq:equi_glob}
\begin{equation*}
\sum_s\passem^{(s)}{\hlam}^{(s)}_b=0
\end{equation*}
%\end{equation}
In order to prove that ${\hlam}\s_b$ also satisfies Fredholm's alternative, we note that since $\kernel\s$ is a basis of $\mathrm{ker}(\stiff\s)$ and $\stiff_{ii}\s$ is invertible, we have $\schur\s\kernel\s_b=\mathbf{0}$ and $\kernel\s_i=-\left.\stiff_{ii}\s\right.^{-1}\stiff_{ib}\s\kernel\s_b$. The condition \label{eq:salambda2} then writes in an equivalent condensed form:
\begin{equation*}
\begin{aligned}
 \left.\kernel\s_b\right.^T \left({\hlam}\s_b+\rhs\s\right)&=0\\
 \left.\kernel\s_b\right.^T \left(\schur\s\dep_b-\rhs\s+ \left.\psassem\s\right.^T\llceil \lam_b \rrceil+\rhs\s\right)&=0\\
% \left(\psassem\s\kernel\s_b\right)^T \left[\left[ \lam_b \right]\right]& =0
\end{aligned}
\end{equation*}
Using the symmetry of $\schur\s$ (inherited from the symmetry of $\stiff\s$) to nullify $\left.\kernel\s_b\right.^T\schur\s$, the condition writes:
\begin{equation}\label{BDD_cond}
\begin{aligned}
\left(\psassem\s\kernel\s_b\right)^T \llceil \lam_b \rrceil=0 
%\, ; \,  \left[\left[ \lam_b \right]\right]
\end{aligned}
\end{equation}
which is exactly the balancing condition \cite{mandel:93:ddm} of the iterative BDD solver: the residual of the BDD iterative solver $\llceil \lam_b \rrceil=\left(\sum_s\passem^{(s)}\lam^{(s)}_b\right)$ \eqref{eq:primal_sys} has to be orthogonal to all local weighted rigid body motions so that preconditioning step is well posed.

Then we have constructed a pair of interface nodal  Vectors $(\hdep_b,\hlam_b)$ which satisfy all required conditions to build admissible fields. 

Note that all the involved operations are already realized during classical steps of the primal domain decomposition approach with a Neumann-Neumann preconditioner and the associated coarse problem, so that all finite element quantities (even the internal ones) are available at no cost; the only extra operations are due to the use of Algorithms $\mathcal{G}_h$ (to compute $\admiss{F}_{bh}$) and $\mathcal{F}_h$ (to compute $\admiss{\sigma}_h$).
% 
% Two remarks need to be made before closing this section:
% \begin{itemize}
% \item The  operations involved in  our method are  exactly the ones  used during
%   preconditioning  the  primal  approach  with  Neumann-Neumann  preconditioner,
%   assuming this excellent preconditioner is  used, our strategy implies no extra
%   cost.
% \item In equation \eqref{eq:bdd_pi} the displacement $\tilde{\dep}\s$ is defined up to a rigid body motion which needs not to be determined since only its symmetric gradient is used to compute the associate stress field $\tilde{\sigma}_h\s$. %the SA stress field.
% \end{itemize}

\section{Recovery of admissible fields in FETI}
\label{sec:fetirecovery}
In the Finite Element Tearing and Interconnecting domain decomposition \cite{FARHAT:1994:ADV}, a unique interface effort unknown $\lam_b$ is introduced so that interface equilibrium is always insured:
\begin{equation}\label{eq:FETI_global}
\lam_b\s={\dassem\s}^T \lam_b \Longrightarrow \sum_s\passem\s\lam_b\s = \mathbf{0}
\end{equation}
Displacements can be deduced from $\lam_b$ if it satisfies Fredholm's alternative on each substructure:
\begin{equation}\label{eq:feti1}
\begin{aligned}
\dep\s &= {\stiff\s}^{+}\left(\force\s+{\traceh\s}^T{\dassem\s}^T\lam_b\right)+\kernel\s\alp\s\\
\mathbf{0}&={\kernel\s}^T\left(\force\s+{\traceh\s}^T{\dassem\s}^T\lam_b\right)
\end{aligned}
\end{equation}
where $\alp\s$ is the unknown magnitude of rigid body motions. 
The FETI solver consists in iteratively finding an interface effort $\lam_b$, under the previous constraint, which insures the  continuity of interface displacement:
%in iteratively trying to verify the continuity of interface displacements under the previous constraints:
\begin{multline}
\mathbf{0}=\sum_s \dassem\s\dep_b\s
= \left(\sum_s \dassem\s\traceh\s{\stiff\s}^+{\traceh\s}^T{\dassem\s}^T\right)\lam_b \\+ \left(\sum_s \dassem\s\traceh\s{\stiff\s}^+\force\s\right)+\left(\sum_s \dassem\s\traceh\s\kernel\s\alp\s\right)
\end{multline}

\subsection{Recovery of SA fields}
In the FETI solver, the nodal interface fields ${\lam}\s_b={\dassem\s}^T\lam_b$ are by construction always balanced at the interface \eqref{eq:FETI_global} and associated to well-posed discrete Neumann problems on each substructure \eqref{eq:feti1}. 
Hence, we can directly set $\hlam_b\s={\lam}\s_b$ and apply algorithms $\mathcal{G}_h$ and $\mathcal{F}_h$ to compute $\admiss{\sigma}\s_h \in \SA{\Omega\s}$ with $\admiss{\sigma}_h=\assemg (\admiss{\sigma}_h\g) \in \SA{\Omega}$.
% \begin{equation*}
% {\sigma}_h\s = \hooke:\varepsilon\left({\dep}_h\s\right) = \hooke:\varepsilon
% \left( \shapev\s{\stiff\s}^+{\traceh\s}^T \left(\force\s+{\traceh\s}^T{\lam}\s_b \right) \right) 
% \end{equation*}
% From which, we are able to compute an admissible stress field $\admiss{\sigma}_h\s$ on each subdomain according to the retained algorithm $\mathcal{F}(.,.)$:
% \begin{equation*}
% \admiss{\sigma}_h\s=\mathcal{F}({\sigma}_h\s, {F}_h\s) \in \SA{\Omega\s, {F}_h\s},
% \end{equation*}
%By construction, the associate assembly stress field is an admissible stress field for the reference problem:
%\begin{equation*}\admiss{\sigma}_h=\assemg (\admiss{\sigma}_h\g) \in \SA{\Omega, g}\end{equation*}

\subsection{Recovery of KA fields}
For a given balanced nodal interface traction $\lam_b$, we introduce, in agreement with \eqref{eq:feti1}, the gap of the interface displacement :
\begin{equation*}
\llfloor \dep_b \rrfloor= \sum_s\dassem^{(s)}\dep^{(s)}_b
\end{equation*}
and we define
\begin{equation*}
\hdep\s_b = \dep\s_b- \left.\dsassem\s\right.^T\llfloor \dep_b \rrfloor
\end{equation*}
where $(\dsassem\s)_s$ are scaled assembling operators so that $\sum_s\dassem\s \left.\dsassem\s\right.^T= \mathbf{I}$. Similarly to the BDD case, a typical  example of such operator $\dsassem\s$  is the multiplicity scaling:
\begin{equation}
  \left.\dsassem\s\right.^T = \left.\dassem\s\right.^T\left(\sum_j\dassem^{(j)}\left.\dassem^{(j)}\right.^T\right)^{-1}
\end{equation}
Note that in  the  case  of  two subdomains, we have:
$\left.\dsassem\s\right.^T \llfloor \dep_b \rrfloor=\displaystyle \frac{1}{2}\llfloor \dep_b \rrfloor$. The connection between FETI and BDD scaling operators (even in the heterogeneous case) is given in \cite{GOSSELET:2003:IEI}.

% In order to compute a displacement field coherent with interface displacement field $\h\dep\s_b$, we consider the problem with prescribed Dirichlet condition on the interface on each subdomain, the solution of which writes:
% \begin{equation}
%   \tilde{\dep}_i\s={\stiff\s_{ii}}^{-1}\left(\force_i\s-{\stiff\s_{ib}}\tilde{\dep}\s_b \right)
% \end{equation}
% \begin{equation}
% \begin{aligned}
%   \tilde{u}_h\s &= \shapev\s\tilde{\dep}\s
%   \end{aligned}
%\end{equation}
It is clear that by construction
\begin{equation*}\sum_s\dassem^{(s)}\hdep\s_b=\mathbf{0}\end{equation*}
Hence nodal interface displacement $\hdep\s_b$ can be used to deduce an admissible displacement field $\admiss{u}\s_h$ so that $\admiss{u}_h=\assemg\left(\admiss{u}_h\g\right)\in\KA{\Omega}$.

Then we have constructed a pair of interface nodal  Vectors $(\hdep_b,\hlam_b)$ which satisfies all required conditions to build admissible fields. Note that all the involved operations are already realized during classical steps of the dual domain decomposition approach (with built-in coarse problem) with Dirichlet's preconditioner, so that all finite element quantities (even the internal ones) are available at no cost: the  quantity   $\llfloor \dep_b \rrfloor$   is  directly available  during the  classical solution  procedure (without  computing any  $\alp^{(j)}$)  which  is   based  on  an  initialization/projection  algorithm  \cite{FARHAT:1994:ADV}, and  the displacement field ${\dep}\s$  can be defined  up  to an  element of  the kernel  (a  rigid body  motion) since  only its  symmetric gradient is used during the computation of the error. The only extra operations are due to the use of algorithms $\mathcal{G}_h$ (to compute $\admiss{F}_{bh}$) and $\mathcal{F}_h$ (to compute $\admiss{\sigma}_h$). 

\section{Numerical assessment}
\label{sec:numeric}
%%%%%%%%%%%%%%%%%%%%%%%%%%%%%%%%%%%%%%%%%%%%%%%%%%%%%%%%%%%%%%%%%%%%%
%%%%%%%%%%%%%%%%%%%%%%%%%%%%%%%%%%%%%%%%%%%%%%%%%%%%%%%%%%%%%%%%%%%%%
%%%%%%%%%%%%%%%%%%%%%%%%%%%%%%%%%%%%%%%%%%%%%%%%%%%%%%%%%%%%%%%%%%%%%

\begin{figure}[ht]
  \centering
\subfigure[Finite element problem ($h=~\frac{L}{4})$\label{fig:gammaStruct}]{
    \includegraphics[height=6.5cm]{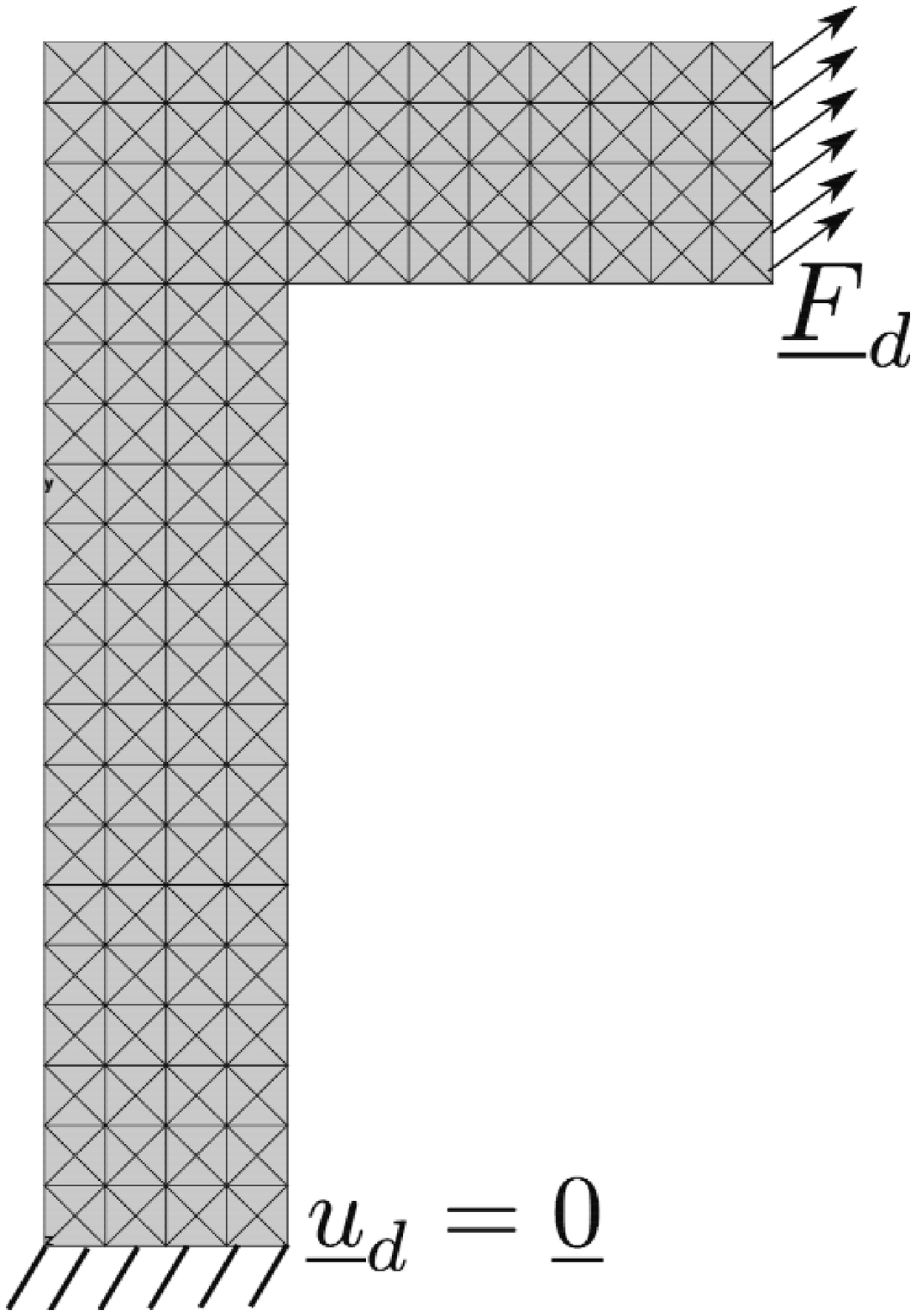}}\qquad
\subfigure[Substructuring ($h=~\frac{L}{8}$, $N_{sd}=8$)\label{fig:gammaStruc8sd}]{
    \includegraphics[height=6.1cm]{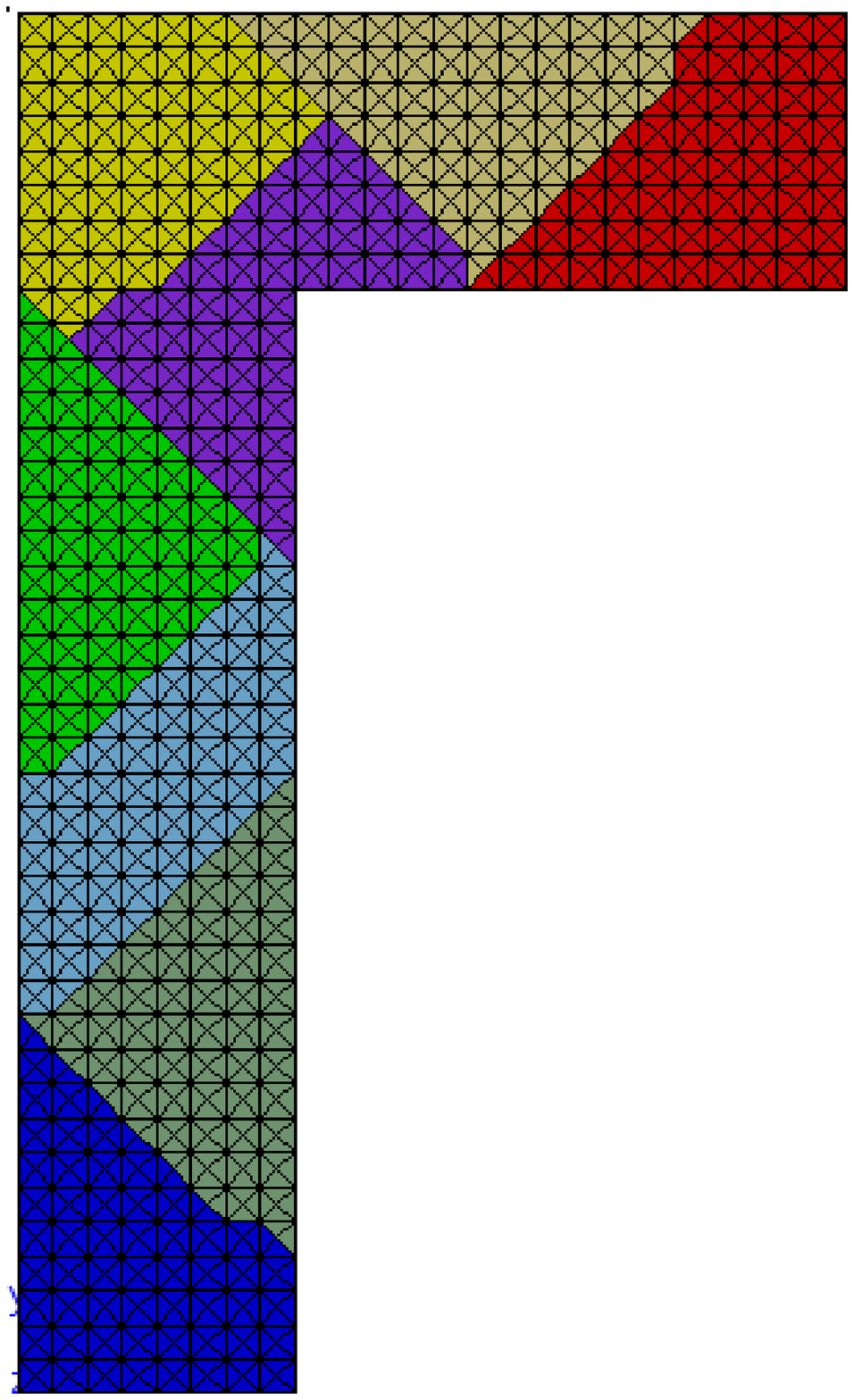}}
  \caption{$\Gamma$-shape structure}
  \label{fig:gammaStructGen}
\end{figure}

In order to assess the performance of our parallel error estimator, we consider the 2D toy problem of the $\Gamma$-shape structure of  Figure \ref{fig:gammaStruct} which has been used in other papers like \cite{chamoin:ladeveze:2008:nonIntrusiveBounds}. Plane stresses are assumed. The material behavior is isotropic, linear and elastic, with Young modulus $E=2000$ MPa and Poisson's ratio $\nu = 0.3$.
The structure is clamped on its basis (whose length is denoted  $L$) and it is submitted to traction and shear on 
its upper-right side, while all the remaining boundaries are traction-free.

Several regular  meshes have  been generated,  constituted  by triangular
elements  of characteristic size  $h=\frac{L}{m}$ with  $m=2,4,8,16,32$. %Let us
%note  that the  figure \ref{fig:gammaStruct}  represent the  mesh  associated to $d=4$. 
For  each mesh, a sequential (mono-domain) computation  is driven, followed  by domain decomposition computations obtained  by an automatic splitting of  the mesh in an increasing number $N_{sd}$ of subdomains ( $N_{sd}=2,4,8$ when $m\leqslant 4$ and
$N_{sd}=2,4,8,16,32$ when $m\geqslant8$).  Figure \ref{fig:gammaStruc8sd} shows such a decomposition for $N_{sd}=8$ and $m=8$.

All the computations   are   driven   in the ZeBuLoN   finite   element   code \cite{ZEBUUSER:2001}, using elements  of polynomial degree $p=1$.
Both  BDD and FETI  algorithms are  used to  solve the  substructured problems,
respectively used  together with Neumann-Neumann and  Dirichlet preconditioners. Beside,
the convergence  criterion of  the solver, which  stands here for  the interface traction gap
 (resp. displacement gap)  in the primal (resp. dual) approach, is set to
$10^{-6}$.

%At  last, for  each computation,  discretization error  is computed  using $p+3$ refinement when  solving the local problems on  elements \cite{babuska:94:vpeena,ladeveze:chamoin:florentin:2009:erdc:sarecovery}.
On each case, in addition to the new parallel error estimator  $\ecrpara$, we compute the standard sequential  $\ecrseq$ and the true error $\globalerror$ obtained using a reference field $u_{ex}$ computed on a very fine mesh:
\begin{equation*}
\begin{aligned}
  \ecrseq &=\ecr{\admiss{u}_h,\admiss{\sigma}_h}{\Omega} \\
  \ecrpara &=\sqrt{\sum_s \left(\ecr{\admiss{u}\s_h,\admiss{\sigma}\s_h}{\Omega\s}\right)^2}\\
\globalerror&= \strainnorm{\strain{u_{ex}-\hat{u}_h}}{\Omega} = \sqrt{{\strainnorm{\strain{u_{ex}}}{\Omega}^2} - {\strainnorm{\strain{\hat{u}_{h}}}{\Omega}^2}}
\end{aligned}
\end{equation*}

%but to  simplify  the comparison, we define  a relative global  discretization error
%$\ecrpara$:
%\begin{equation}
%  \label{eq:erel}
%  \ecrpara=\frac{\ecr{\admiss{u}_h,\admiss{\sigma}_h}{\Omega}}{\strainnorm{\varepsilon{(u_{ex})}}{\Omega}}
%\end{equation}
%where $u_{ex}$ stands for the exact solution of the problem, which has been computed separately using a very fine mesh.

\subsection{Quality of the parallel error estimator}
\begin{figure}[ht]
  \centering
% \begin{tabular}[h]{cc}
%  \psfrag{relative error}{$\ecrpara$}
%  \includegraphics[width=0.45\textwidth]{figures/H_conv/D_error_h} & 
% \psfrag{relative error}{$\ecrpara$}
%  \includegraphics[width=0.45\textwidth]{figures/H_conv/P_error_h}\\
%FETI Case  &   BDD Case
%&
% \psfrag{relative error}{$\ecrpara$}
%\psfrag{Error vs. discretization}{}
\includegraphics[width=0.75\textwidth]{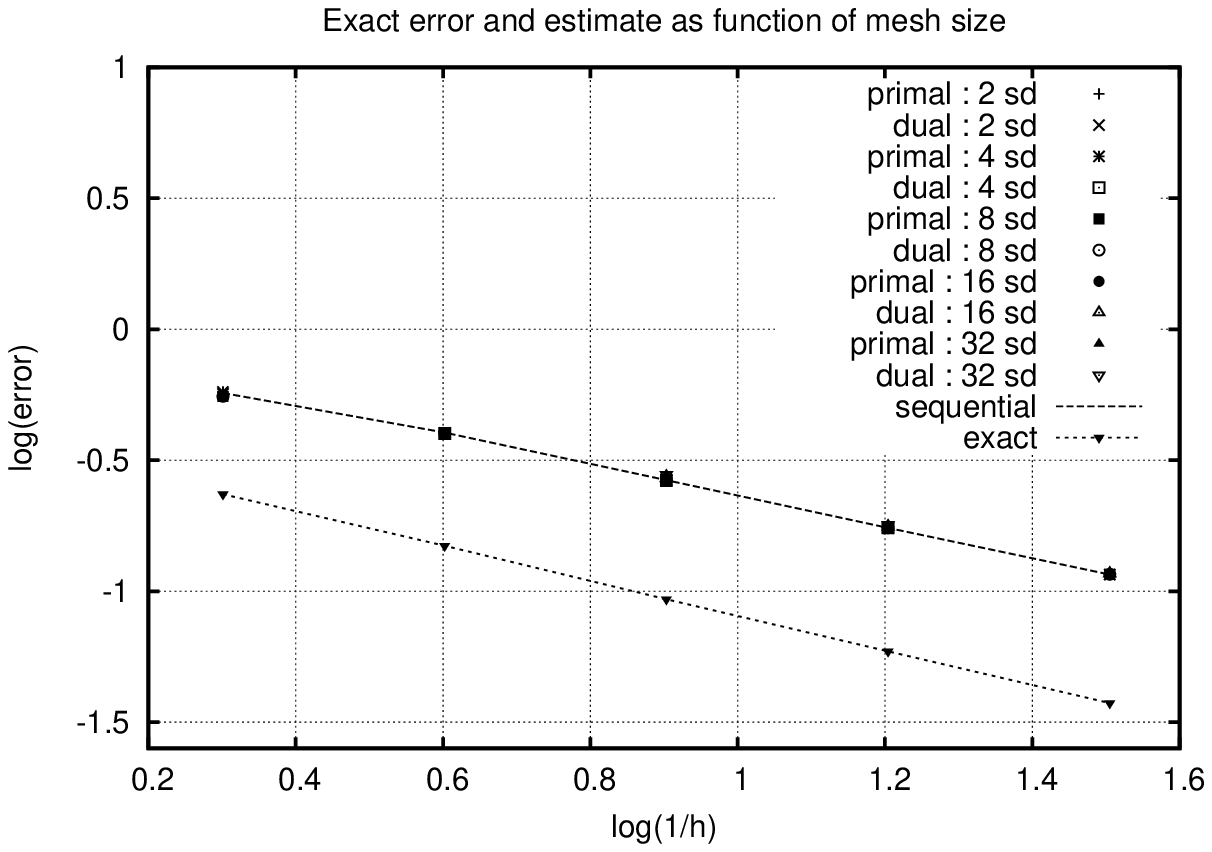}
%\end{tabular}
\caption{Convergence of Error $\globalerror$ and Estimators $\ecrseq$ and $\ecrpara$ (for various $N_{sd}$) vs. element size $h$ }
  \label{fig:error_vs_h}
\end{figure}

We first study the quality of the parallel error estimator $\ecrpara$ for computations when convergence of the  domain decomposition solver is reached. As said earlier, the proposed technique does not lead to the same statically admissible field because of the special treatment of the interface traction \eqref{eq:Fh}. Our estimator might then be sensitive to the substructuring, we thus compare the estimations obtained with meshes of characteristic size $h$ and decomposition into $N_{sd}$ subdomains. Results are given in Figure \ref{fig:error_vs_h} and Table \ref{tab:erel}.

\begin{table}[ht]
  \centering
  \begin{tabular}[ht]{cccccccc}
    \hline
    $h$ & $L/2$ & $L/4$ & $L/8$ & $L/16$ & $L/32$ \\
    \hline
    $\#$ dofs & 146 & 514 & 1922 & 7426 & 29186 \\
    \hline
    \hline
    $\globalerror$ & 0.2347 & 0.1493 & 0.0937 & 0.0597 & 0.0386 \\
    \hline
    $\ecrseq$ & 0.5712 & 0.4035 & 0.2662 & 0.1769 & 0.1151 \\
    \hline
    \hline
    $N_{sd}$ & \multicolumn{5}{c}{$\ecrpara$} \\
    \hline
2	&0.5657	&0.4021	&0.2648	&0.1747	&0.1151\\
4	&0.5768	&0.4007	&0.2648	&0.1747	&0.1151\\
8	&0.5546	&0.4007	&0.2676	&0.1747	&0.1165\\
16	&	&	&0.2690	&0.1761	&0.1165\\
32	&	&	&0.2787	&0.1789	&0.1178\\
    \hline
  \end{tabular}
  \caption{Error $\globalerror$ and Estimators $\ecrseq$ and $\ecrpara$ (for various $N_{sd}$) vs. element size $h$}
  \label{tab:erel}
\end{table}

%\begin{table}[h!t]
%  \centering
%  \begin{tabular}[ht]{cccccccc}
%    \hline
%    $h_e$ & $L/2$ & $L/4$ & $L/8$ & $L/16$ & $L/32$ \\
%    \hline
%    $\#$ dofs & 146 & 514 & 1922 & 7426 & 29186 \\
%    \hline
%    \hline
%    $N_{sd}$ & \multicolumn{5}{c}{$\ecrpara$} \\
%    \hline
%    1 & 2.52 & 1.52 & 0.93 & 0.59 & 0.38 \\
%    \hline
%    2 & 2.50 & 1.52 & & 0.59 & 0.38 \\
%    4 & 2.50 & 1.52 & 0.93 & 0.59 & 0.38 \\
%    8 & 2.45 & 1.52 & 0.94 & 0.59 & \\
%    16 & & & 0.95 & 0.59 & 0.38 \\
%    32 & & & 1.00 & 0.61 & \\
%    \hline
%  \end{tabular}
%  \caption{Behavior of $\ecrpara$ as a function of $h_e$ and $N_{sd}$ (FETI solver)}
%  \label{tab:erel}
%\end{table}

We observe that:
\begin{itemize}
  \item The results obtained by FETI and BDD can not be distinguished (which is why only FETI results are given in Table \ref{tab:erel}).
  \item $\ecrpara$ barely depends on the substructuring; the results are quite similar whether they are conducted on a single domain (``sequential'' curve) or on $N_{sd}$ subdomains. Only a slight rise of the estimation can be observed when the number of interface degrees of freedom is not small compared to the number of internal degrees of freedom, which is logical since the description of interface traction fields is coarser in parallel than in sequential.
\end{itemize}

As a conclusion, the parallel error estimator $\ecrpara$ enables to recover the same efficiency factor as the standard sequential one, while the CPU-time is divided by $N_{sd}$.

\subsection{Convergence of the parallel estimator along DD-solver iterations}
Previous results enabled to analyse the quality of the parallel estimator when interface quantities had converged. A new feature associated to the use of an iterative solver for the domain decomposition (DD) problem is that the discretization error estimation can be conducted before DD convergence is reached, that is in presence of  displacement or traction  discontinuity at the interface as explained in Sections~\ref{sec:bddrecovery} and \ref{sec:fetirecovery}

\begin{figure}[ht]
\centering
\begin{tabular}[h]{cc}
  \psfrag{error estimator}{$\ecrpara$}
  \includegraphics[width=0.465\textwidth]{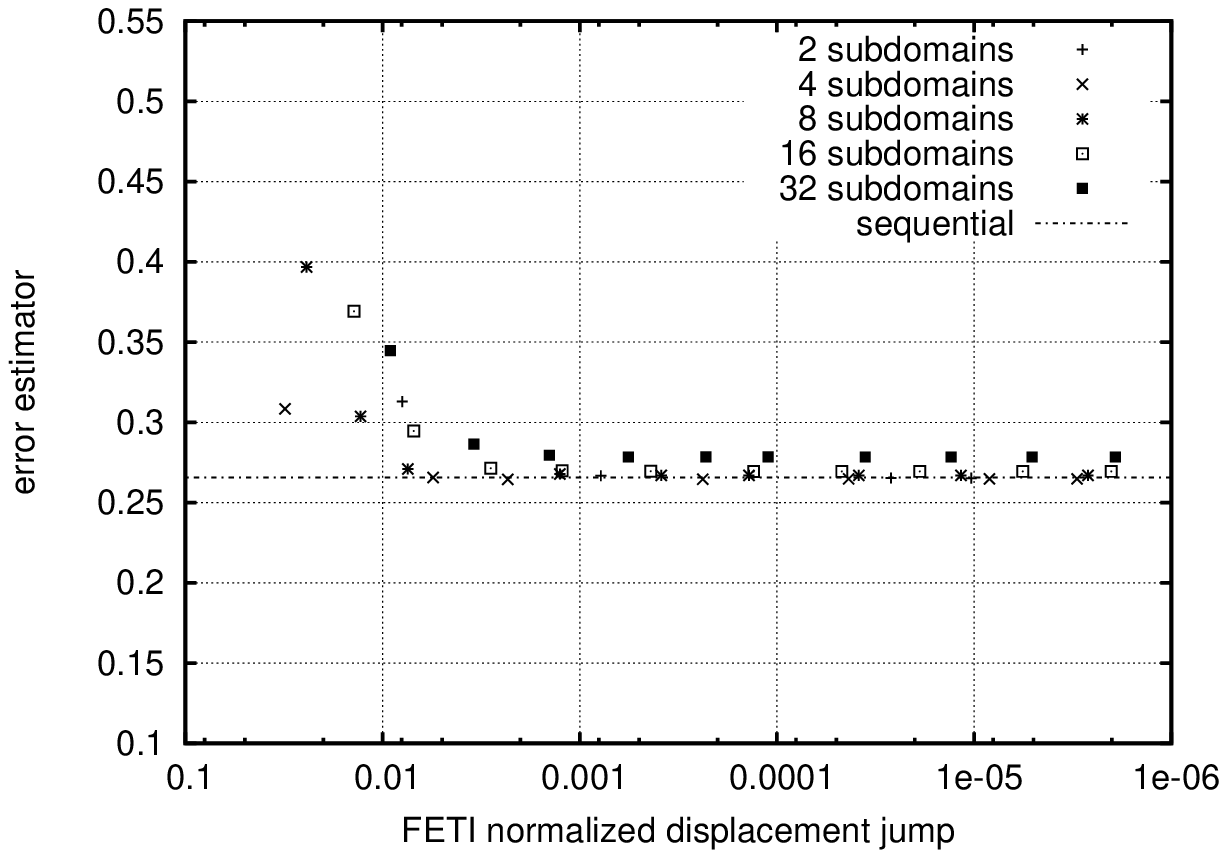} 
 &\psfrag{error estimator}{$\ecrpara$}
  \includegraphics[width=0.465\textwidth]{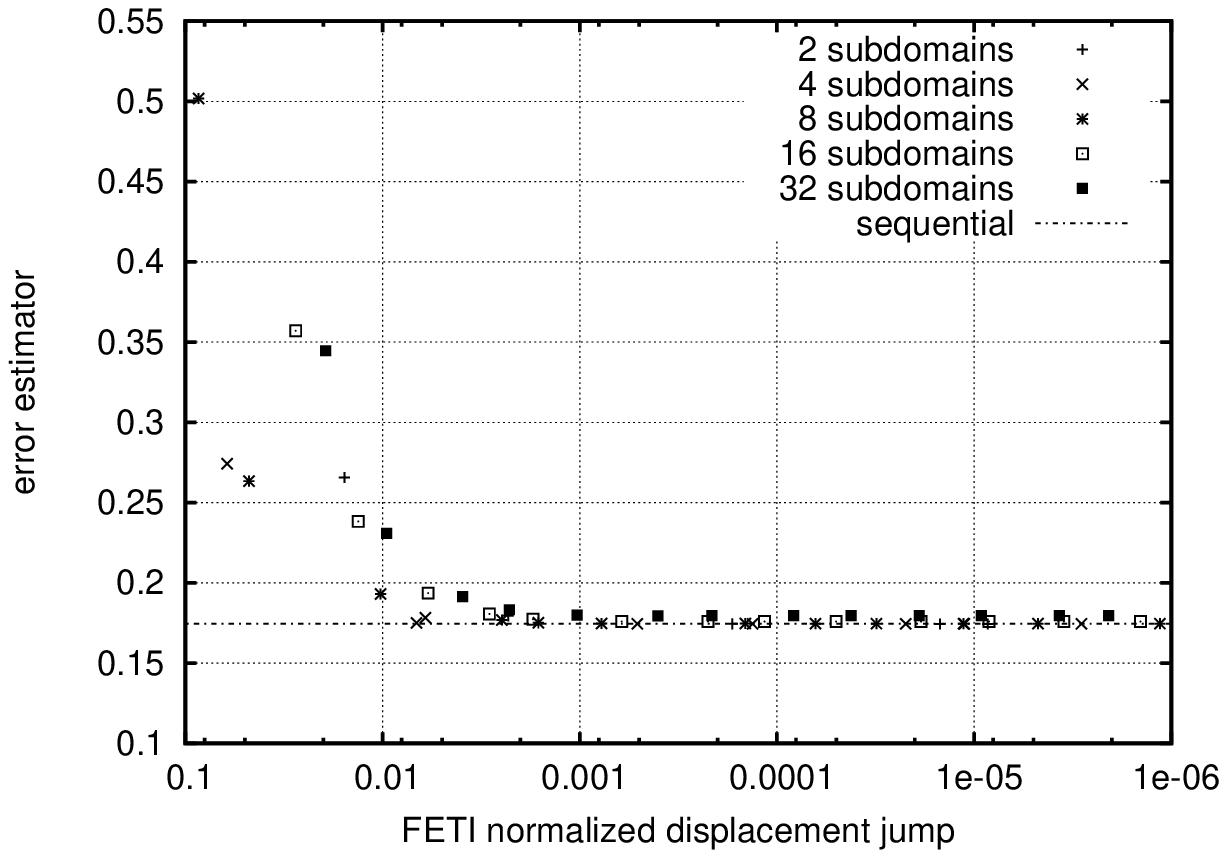}\\
FETI Case: $h=L/8$  & FETI Case: $h=L/16$\\
& \\
  \psfrag{error estimator}{$\ecrpara$}
  \includegraphics[width=0.465\textwidth]{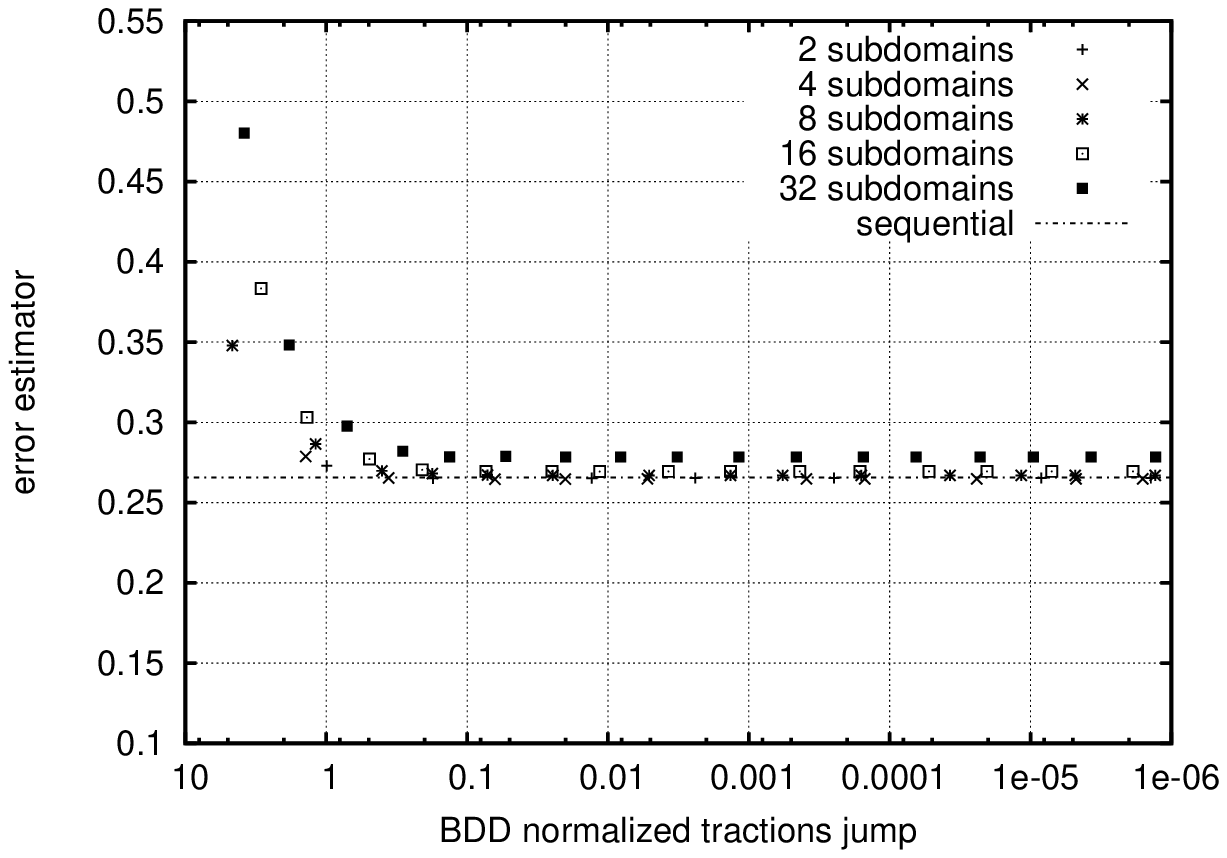}
 &\psfrag{error estimator}{$\ecrpara$}
  \includegraphics[width=0.4650\textwidth]{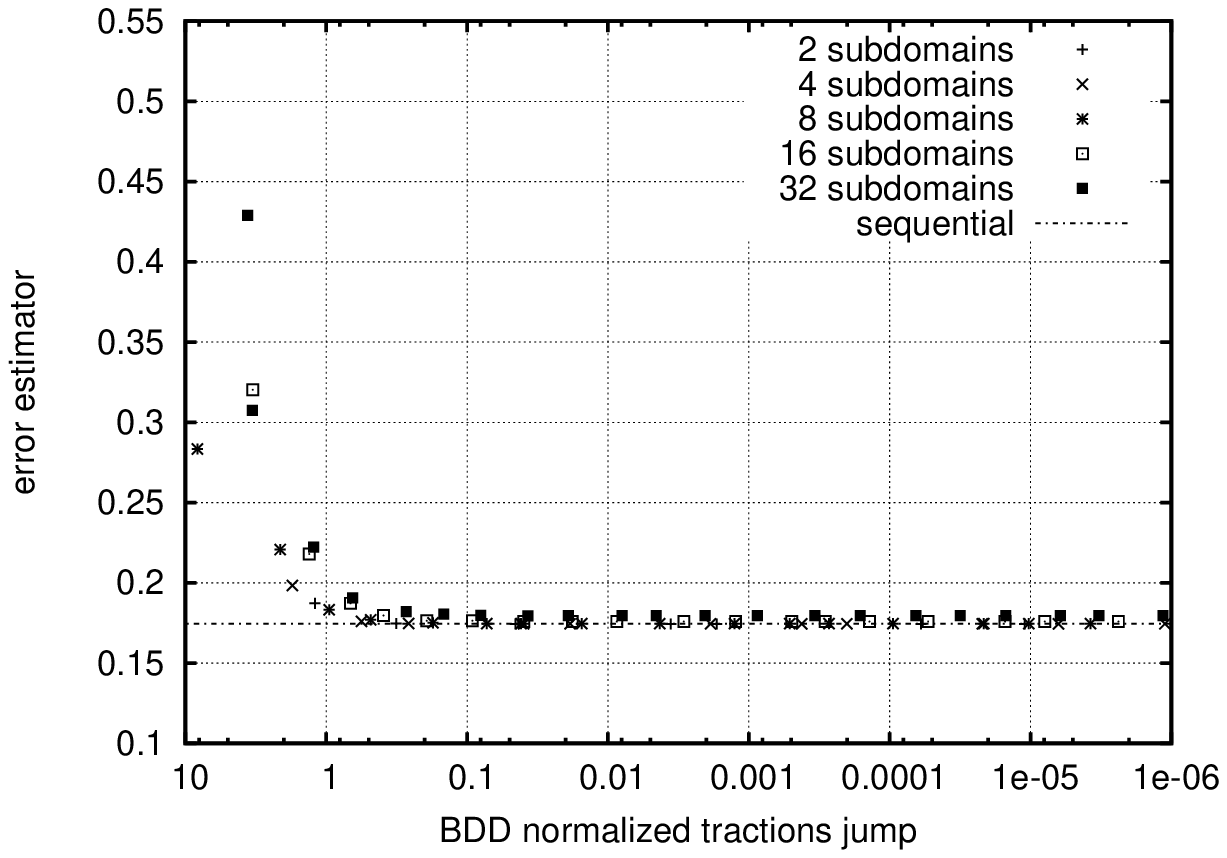}\\
BDD Case: $h=L/8$  & BDD Case: $h=L/16$\\
%& \\
%  \psfrag{relative error}{$\ecrpara$}
%  \includegraphics[width=0.45\textwidth]{figures/Erreur_Residu/DP_error_gap_h08}
% &  \psfrag{relative error}{$\ecrpara$}
%  \includegraphics[width=0.45\textwidth]{figures/Erreur_Residu/DP_error_gap_h16}\\
%$h_e=L/8$  & $h_e=L/16$\\
\end{tabular}
\caption{Convergence  of estimator vs. DD residual}% at $h_e=L/8$}
  \label{fig:error_convergence}
\end{figure}

We then compute the parallel error estimator $\ecrpara$ at each iteration of the DD solver. Convergence curves of $\ecrpara$ during the FETI and BDD iterations are shown on Figure \ref{fig:error_convergence}.
Parallel error estimator is plotted as a function of the FETI (resp. BDD) residual, defined (for Iteration $n$) as the normalized displacement (resp. traction) gap at the interface:
\begin{equation}
  r^n=\frac{\llfloor \dep_b^n \rrfloor_{\Gamma}}{\llfloor \dep_b^0 \rrfloor_{\Gamma}} \qquad \text{ or }\qquad 
r^n=\frac{\llceil \lam_b^n \rrceil_{\Gamma}}{\llceil \lam_b^0 \rrceil_{\Gamma}}
\end{equation}
Classical stopping criterion for the of convergence of DD solver is this residual being below $10^{-6}$. Because of the similarity between the curves, the only shown cases correspond to $h=L/8$ and $h=L/16$.

%\begin{figure}[ht]
%  \centering
%  \psfrag{relative error}{$\ecrpara$}
%  \includegraphics[width=0.75\textwidth]{figures/Erreur_Residu/DP_error_gap_h08}
%  \caption{Convergence  of global  error vs. residual at  $h_e=L/8$}
%  \label{fig:error_convergence_h08}
%\end{figure}
%\begin{figure}[ht]
%  \centering
%  \psfrag{relative error}{$\ecrpara$}
%  \includegraphics[width=0.75\textwidth]{figures/Erreur_Residu/DP_error_gap_h16}
%  \caption{Convergence  of global  error vs.  residual at  $h_e=L/16$}
%  \label{fig:error_convergence_h16}
%\end{figure}

The curves show a rapid convergence  of the parallel error estimator along iterations of the  solver,  so that  $\ecrpara$  can be  considered  as  converged when  FETI residual reaches an order of  magnitude of $5.10^{-3}$ or BDD residual reaches $5.10^{-1}$, which corresponds to at most 5 iterations whereas the solver convergence is achieved in 10 to 20 iterations.

Actually, $\ecrpara$ is driven by both the discretization error and the convergence of the solver (interface error). The ``L''-shaped curves show that the impact of residual of the DD solver is preponderant only at the first iterations (when interface fields are very poorly estimated), after $\ecrpara$ stagnates at a value very close to $\ecrseq$ which is only associated to the discretization error.

%The impact of the residual of the solver is preponderant before $\ecrpara$ has converged (at that time one can see that the error is overestimated on the interfaces), whereas the  estimate  only depends on discretization error when the solver residual reaches a  certain order of magnitude (at that time discretization has no influence).

Then, it seems possible to stop the iterations  of the solver far before convergence while still obtaining an accurate global estimate for the discretization error.

\begin{figure}[ht]
  \centering
  \begin{tabular}[h]{ccc}
    \includegraphics[width=2.9cm]{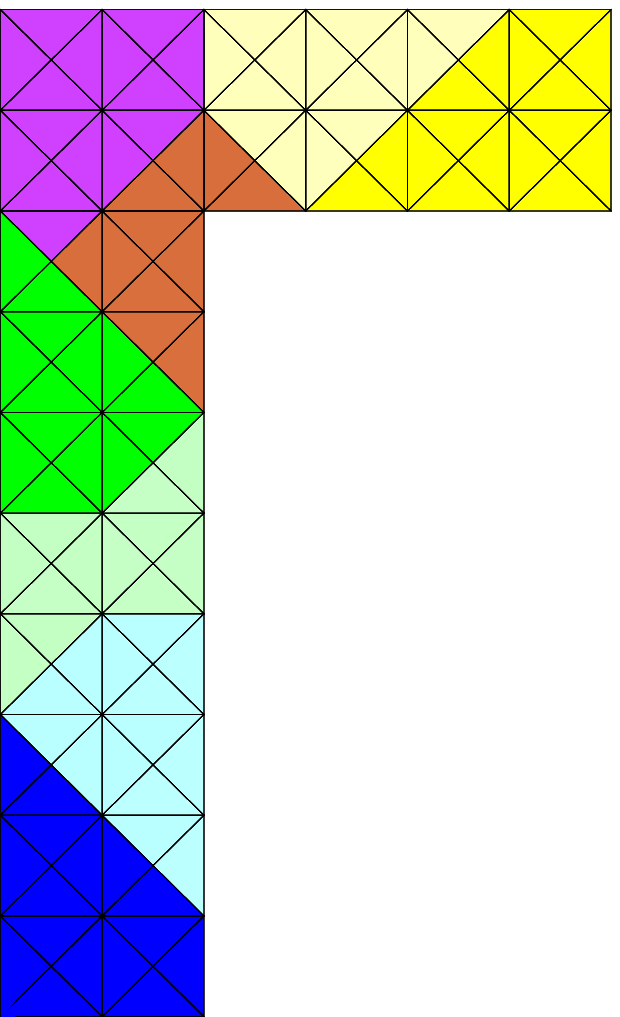}  
    & \includegraphics[width=2.9cm]{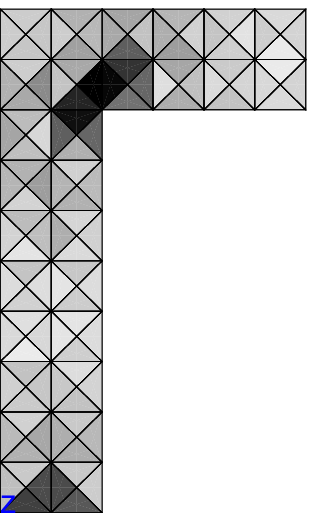}
    & \includegraphics[width=1.5cm]{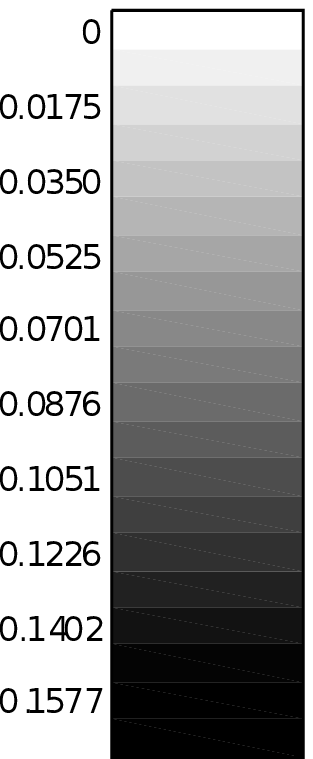}  \\
    Decomposition & Reference map & Range 
    \vspace{1pt} \\
    \includegraphics[width=2.9cm]{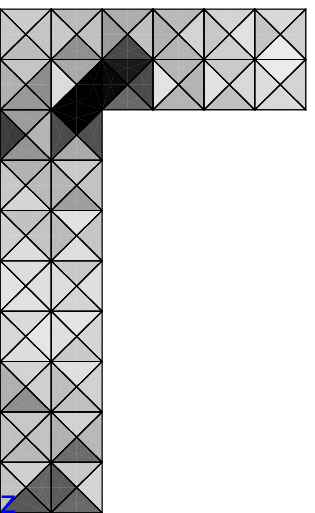}
   & \includegraphics[width=2.9cm]{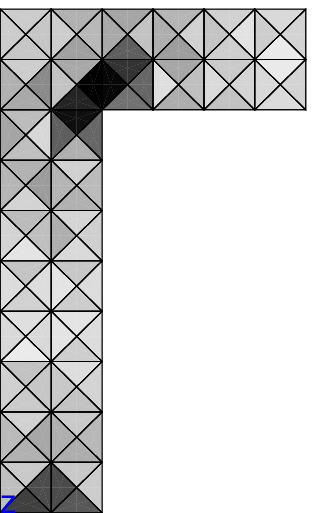}
   & \includegraphics[width=2.9cm]{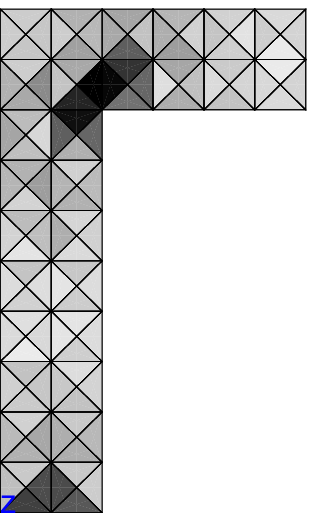}\\
    Iteration 1 &  Iteration 4 & Iteration 5 \\
		&  \textbf{BDD solver} & \\
    \vspace{1pt} \\
     \includegraphics[width=2.9cm]{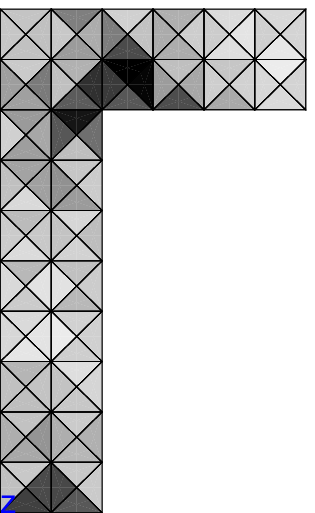}
    & \includegraphics[width=2.9cm]{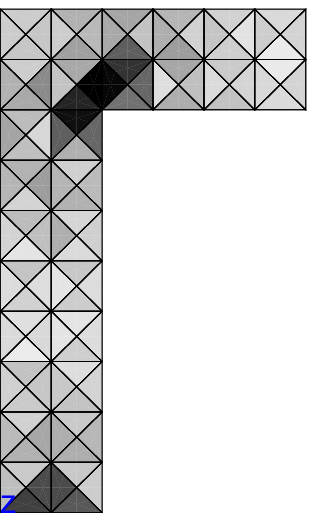} 
    & \includegraphics[width=2.9cm]{h02sd8iter5primal-nb.eps}\\
    Iteration 1 &  Iteration 4 & Iteration 5 \\
		&  \textbf{FETI solver} &
  \end{tabular}
  \caption{Maps of $\ecrparaE$ for $h=L/2$ and $N_{sd}=8$ at various iterations}
  \label{fig:eltmaph2}
\end{figure}

\begin{figure}[ht]
  \centering
  \begin{tabular}[h]{ccc}
   \includegraphics[width=2.5cm]{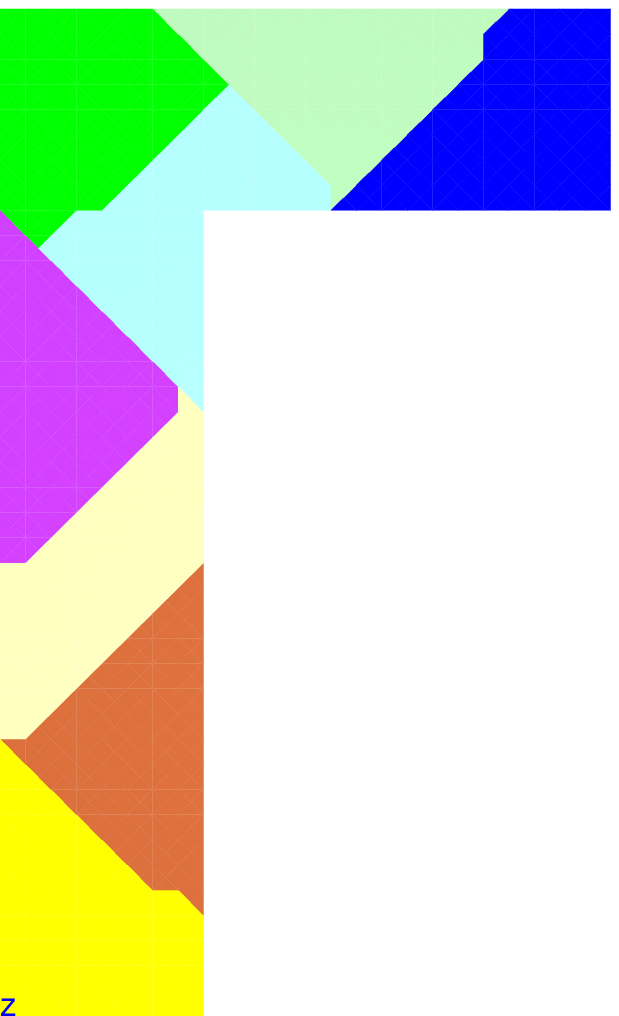}  
    & \includegraphics[width=2.5cm]{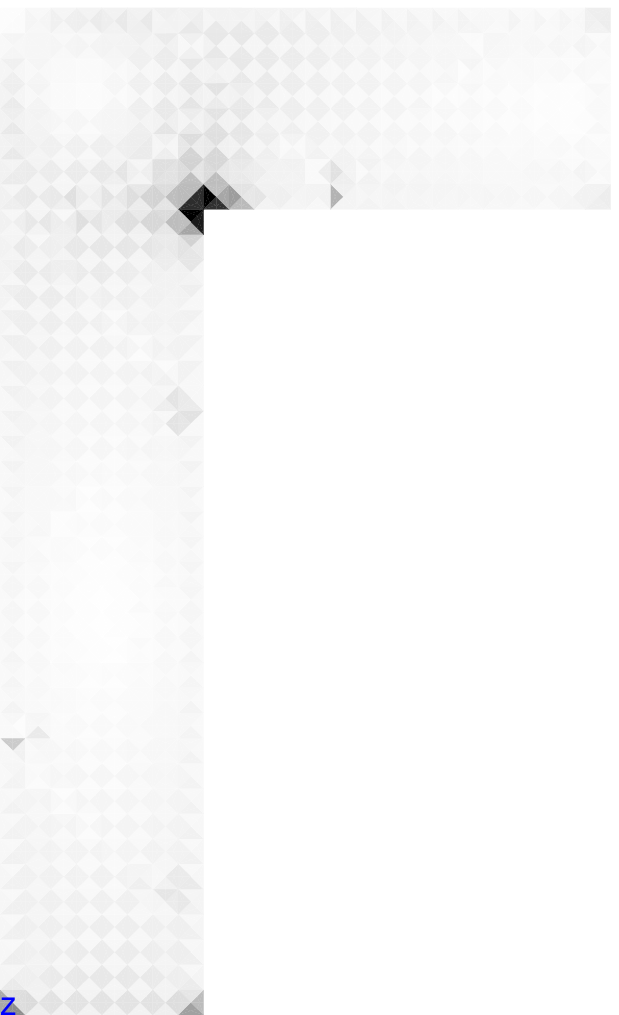}
    & \includegraphics[width=1.5cm]{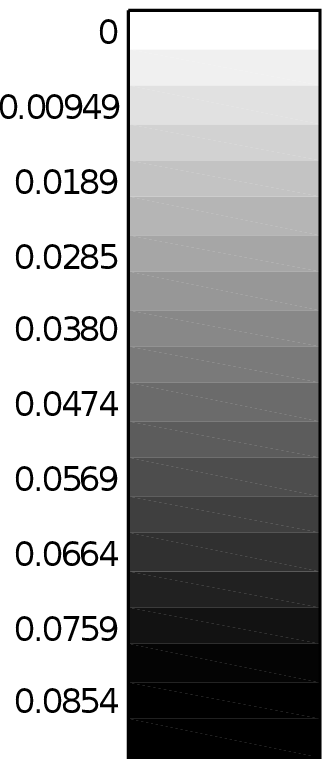}
 \\
        Decomposition & Reference map & Range 
    \vspace{2pt} \\
    \includegraphics[width=3.cm]{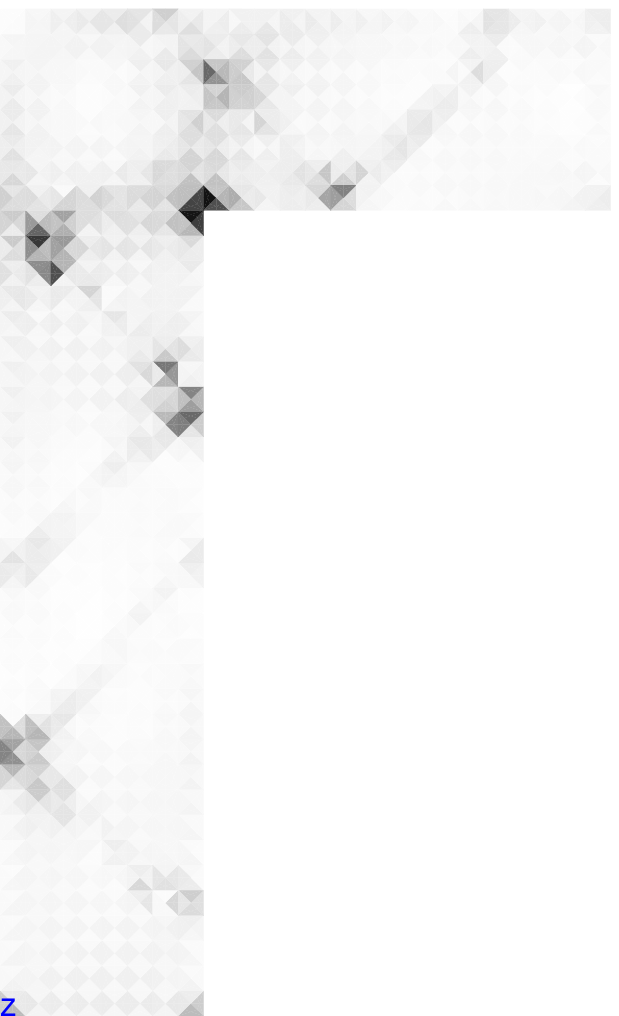}
   & \includegraphics[width=3.cm]{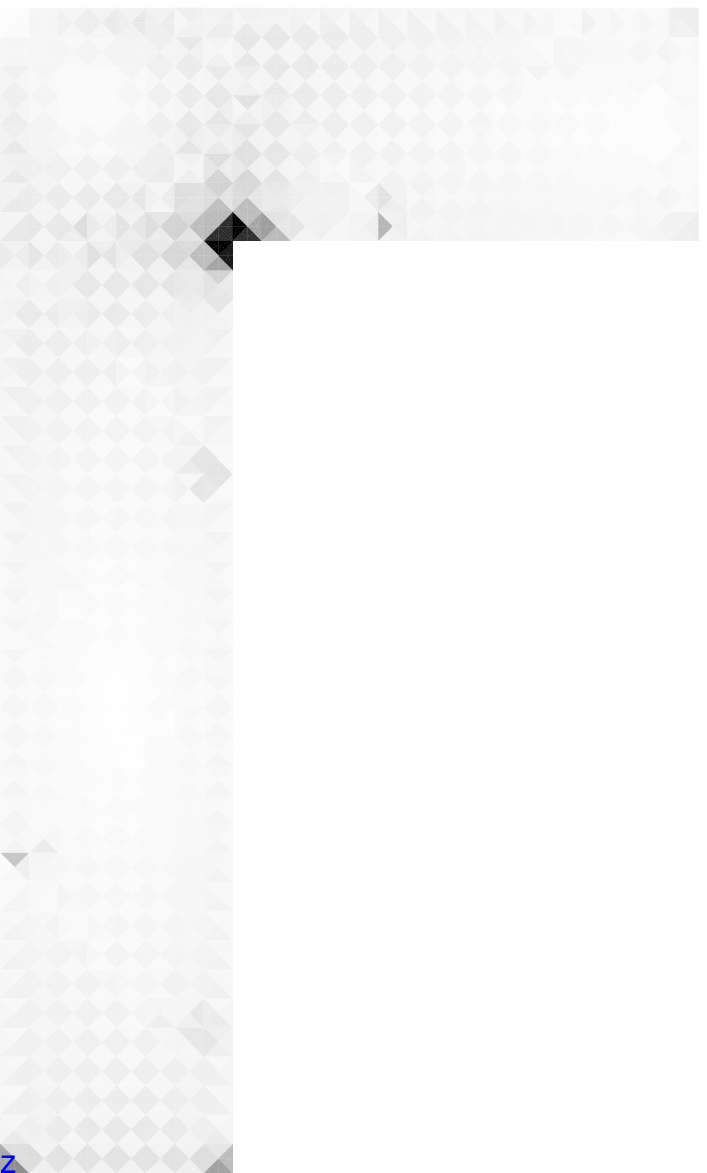}
   & \includegraphics[width=3.cm]{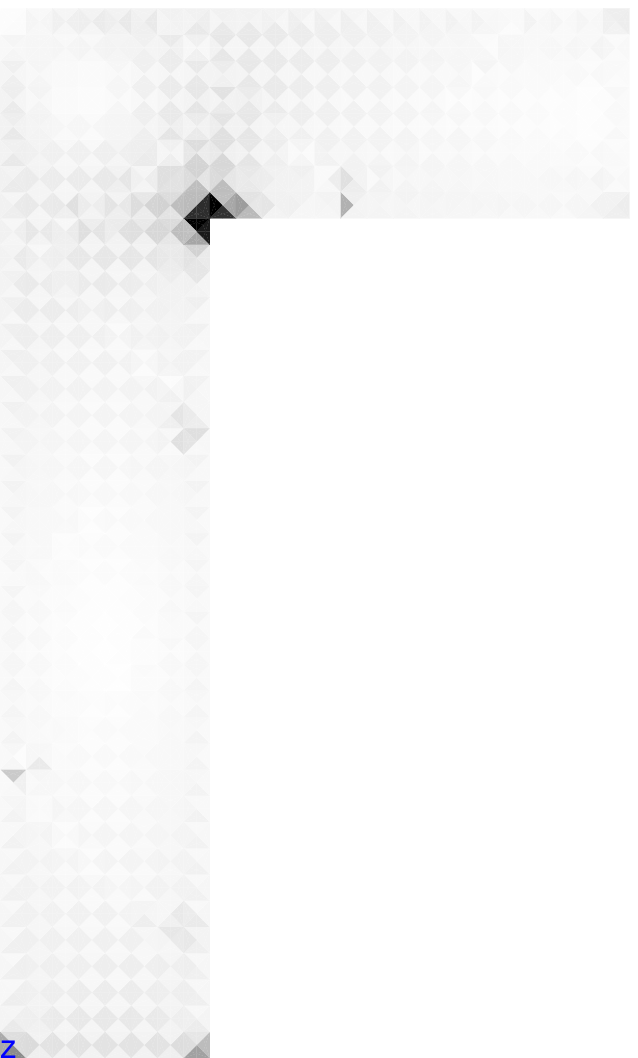}\\
    Iteration 1 &  Iteration 4 & Iteration 5 \\
		&  \textbf{BDD solver} & \\
    \vspace{2pt} \\
     \includegraphics[width=3.cm]{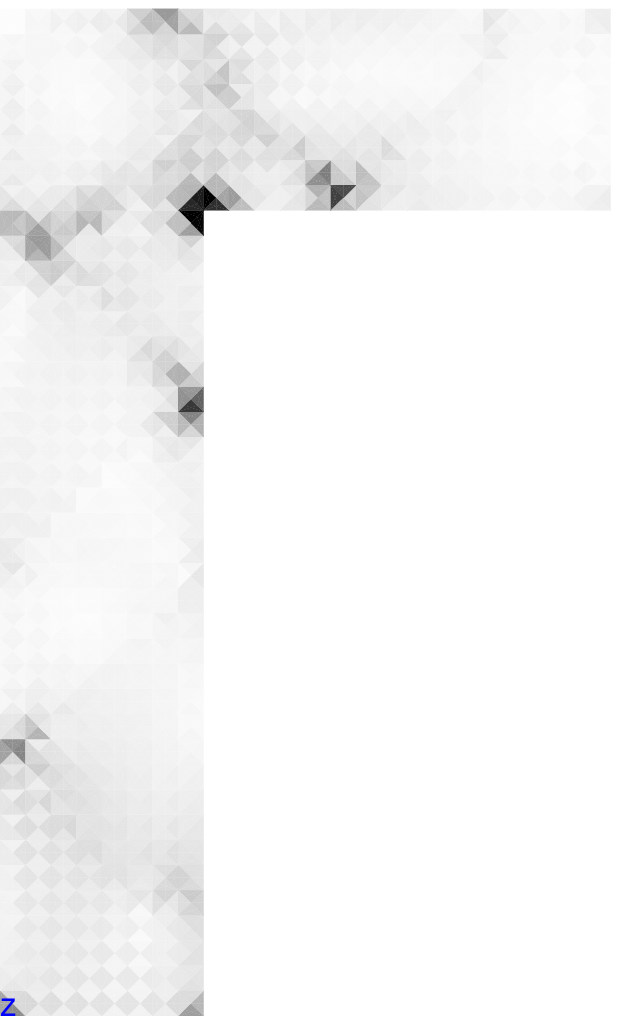}
    & \includegraphics[width=3.cm]{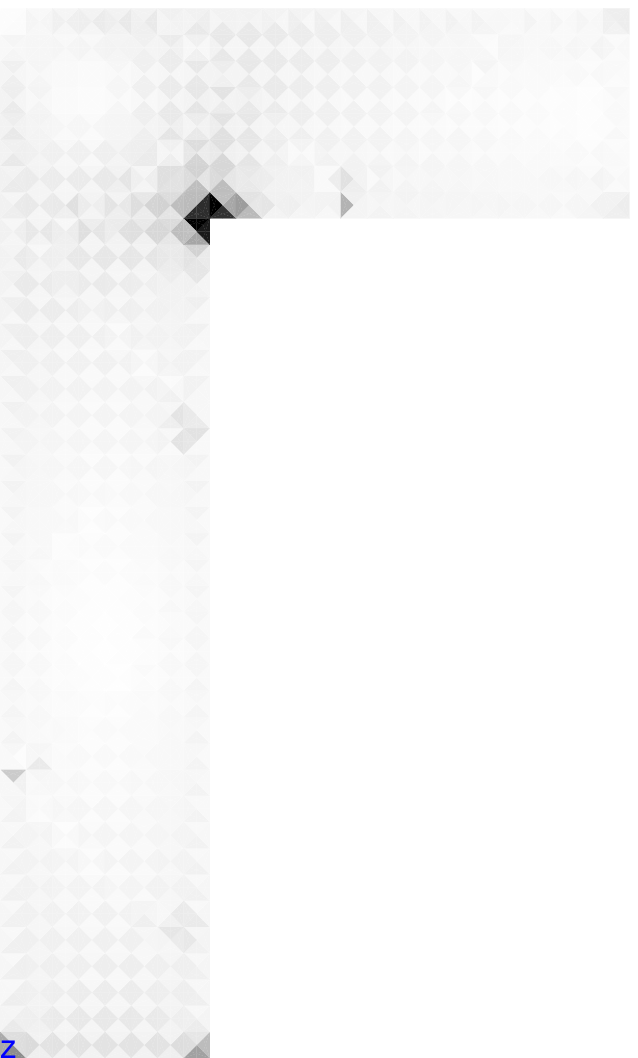} 
    & \includegraphics[width=3.cm]{h08sd8iter5primal-nb.eps}\\
    Iteration 1 &  Iteration 4 & Iteration 5 \\
		&  \textbf{FETI solver} & \\
    \vspace{2pt} \\
  \end{tabular}
  \caption{Maps of $\ecrparaE$ for $h=L/8$ and $N_{sd}=8$ at various iterations}
 \label{fig:eltmaph8}
 \end{figure}

Figures  \ref{fig:eltmaph2} and  \ref{fig:eltmaph8} show maps of  the  elementary contributions  $\ecrparaE$ to the parallel error estimator $\ecrpara$ at
different steps of the convergence, for $N_{sd}=8$ with $h_{e}=L/2$ or $h_{e}=L/8$. At the first iterations it can be seen that the estimator highlights both discretization errors (around the re-entrant angle) and lack of convergence of the solver (along the interfaces), whereas very quickly the solver (that is the interfaces) does not contribute any more to the estimator.

The various examples show that the convergence of the global estimator is due to the convergence of elementary Contributions $\ecrparaE$, which means that when willing to carry out remeshing procedures, the maps obtained after few iterations of the solver are sufficient to define correct refinement instructions.

%Figure show the  same maps  for   and $N_{sd}=8$, with the splitting of the fig. \ref{fig:gammaStruc8sd}.

\section{Conclusions}
\label{sec:conclusions}

In this paper, we presented a new approach to  handle robust model verification based on constitutive relation error in a domain decomposition context. 

The method relies on the construction of fields that are kinematically and statically admissible on the whole structure. We showed that a fully parallel construction is possible even when starting from fields which do not satisfy interface conditions. The construction is a three-step procedure: first displacement and traction nodal fields are built-up so that discrete admissibility conditions are satisfied, second continuous admissible traction fields are deduced, third these fields are used as input by any classical recovery procedure. The first step is implicitly done when good preconditioners are employed within the domain decomposition methods and the second step corresponds to the inversion of small and sparse ``mass'' matrices.

Our first results show that not only the estimation error does not suffer from the approximation that are made at the interface in order to achieve full parallelism, but that even roughly estimated interface fields enable to obtain a good estimation of the discretization error and correct maps of elementary contributions which are required by mesh adaptation procedures. Thus not only the computational cost are divided by the number of processors but the prior obtainment of the finite element solution can be accelerated since a coarse solution is sufficient (which corresponded to 3 to 5 times less iterations in our case).

Future studies will deal with parallel mesh adaptation.

\bibliography{apf}

\end{document}